\documentclass[reqno]{elsart}
\usepackage{graphicx}
\usepackage{subfigure}
\usepackage{amsmath,amssymb}
\usepackage{color}
\definecolor{oneblue}{rgb}{0,0.0,0.75}

\usepackage{dsfont}
\usepackage{textgreek}
\makeatletter
\newcommand{\sech}{\mathop{\operator@font sech}}
\newcommand{\sign}{\mathop{\operator@font sign}}
\makeatother

\newtheorem{lemma}{Lemma}[section]

\newtheorem{remark}{Remark}[section]
\numberwithin{equation}{section}

\renewcommand{\nu}{\text{\textnu}}

\renewcommand{\eta}{\text{\texteta}}
\renewcommand{\beta}{\text{\textbeta}}
\renewcommand{\mu}{\text{\textmugreek}}
\renewcommand{\alpha}{\text{\textalpha}}
\renewcommand{\kappa}{\text{\textkappa}}
\renewcommand{\omega}{\text{\textomega}}
\renewcommand{\theta}{\text{\texttheta}}

\renewcommand{\geq}{\geqslant}

\renewcommand{\sim}{\thicksim}

\begin{document}
\begin{frontmatter}
\title{An asymptotic model for internal capillary-gravity waves in deep water}


\author{A. Dur\'an}
\address{ Applied Mathematics Department,  University of
Valladolid, P/ Bel\'en 15, 47011 Valladolid, Spain. Email:angel@mac.uva.es}


\begin{abstract}
Considered in this paper is a bi-directional model for the propagation of interfacial capillary-gravity waves in a two-layer system of fluids with rigid lid condition for the upper layer and  lower layer with a much larger or infinite depth. The system is derived from a reformulation of the Euler equations for internal waves with nonnegligible surface tension effects in the interface and the corresponding asymptotic model under the Benjamin-Ono regime. Another unidirectional model, so-called regularized Benjamin equation, generalizing the Benjamin equation, is also introduced. Well-posedness of the new equations and existence of solitary wave solutions are discussed.
\end{abstract}
\end{frontmatter}

\section{Introduction}
The present paper is concerned with the derivation and analysis of an asymptotic model for internal capillary-gravity waves. The model incorporates bi-directionality and the physical regime under which the corresponding differential system is derived is compatible with that of established by Benjamin, \cite{Benjamin1967,Benjamin1992,Benjamin1996}, rigurously developed by Albert et al., \cite{ABR}, for the Benjamin equation.

\begin{figure}[htbp]
\centering
{\includegraphics[width=10cm,height=4cm]{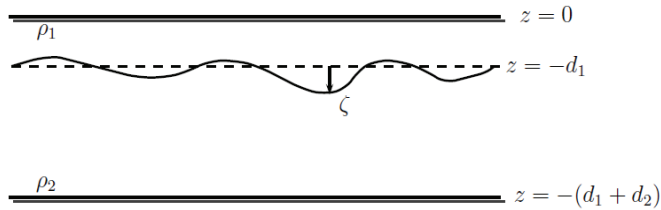}}
\caption{Idealized model of internal wave propagation in a two-layer interface.}
\label{f0}
\end{figure}

More specifically, the idealized model under study consists of two inviscid, homogeneous, irrotational incompressible fluids of depths $d_{1}<<d_{2}$ and different (constant) densities $\rho_{1}<\rho_{2}$, see Figure \ref{f0}. The upper layer satisfies a rigid lid assumption (that is, it is bounded above by an impenetrable, bounding surface) while the lower layer is bounded below by an impermeable, horizontal and flat bottom. The interest is in the description of the motion of the deviation of the interface between the fluids, which is affected by both gravity and capillary forces.

As mentioned by Kalisch in \cite{Kalisch2007}, much of the literature on capillary-gravity interfacial waves with rigid lid for the upper layer and infinitely deep lower fluid concerns analytical or computational studies of solitary waves, from the ones with oscillatory decay admitted as solutions of the weakly nonlinear model derived by Benjamin, \cite{Benjamin1992}. The stability of these waves under small perturbations was computationally analyzed by Calvo and Akylas, \cite{CA}, from the full Euler equations. On the other hand, the computations of interfacial capillary-gravity waves by Laget and Dias, \cite{LagetD1997}, are also based on the numerical approximation of an integro-differential formulation of the full Euler equations (see also \cite{DiasMV1996} and the analytical study by Dias and Iooss, \cite{DiasI1996}). The numerical results are compared with the experimental results by Koop and Butler, \cite{KoopB1981}. Concerning the formulation of asymptotic models, an extension of the Benjamin equation which allows weak transverse variations is derived by Kim and Akylas, \cite{KimA2006}. This is used to study gravity-capillary lumps from the previous numerical results obtained in the first part of the work, \cite{KimA2005}. Finally, Kalisch, \cite{Kalisch2007}, proposes some one-dimensional systems for the propagation of interfacial waves subject to capillary forces. His derivation is based on formal asymptotic expansions of the velocity potential associated to the layers form the one-dimensional Euler equations and the combination with the hypotheses of the physical regime for the Benjamin equation described in \cite{Benjamin1992,ABR}. Indeed this list of references is far from being exhaustive and can be additionally extended if the assumptions for the layers (concerning either the boundary conditions or the depths) are modified, \cite{HelfrichM}.

In a different setting from that used in \cite{KimA2006} and \cite{Kalisch2007}, for the derivation of the model proposed in the present paper the approach developed by Bona et al. in \cite{BLS2008} is considered. This is based on several steps: first, a reformulation of the Euler system for internal waves is made by using two nonlocal operators linking the velocity potentials for the layers at the interface. Then suitable asymptotic expansions of these operators, in accordance with the physical regime under study for the layers, allows to derive the corresponding asymptotic model from the Euler system in some consistent, precisely defined way.

The steps of this approach is adapted in the present paper as follows. We must first consider the Euler system for internal waves which includes surface tension effects at the interface. In \cite{Lannes}, the derivation of the corresponding equations makes use of the Dirichlet-to-Neumann operators associated to the two fluid layers, leading to a system of two equations for the deviation of the interface and a suitable combination of the traces of the velocity potentials at the interface. As mentioned by Lannes, when $\rho_{1}\neq 0$ these are the canonical variables of the Hamiltonian formulation made by Benjamin and Bridges, \cite{BenjaminB}, while if $\rho_{1}=0$, the formulation reduces to tha of the case of surface waves due to Zakharov, \cite{Zakharov}, and Craig and Sulem, \cite{CraigS}. In the present paper, the derivation will be made by using the nonlocal operators considered in \cite{BLS2008}.

Then the formulation of the Benjamin asymptotic model from the resulting Euler system takes into account the physical regime that Benjamin assumed to obtain his uni-directional equation. The validity of this regime was specified by Albert et al. in \cite{ABR}. Their analysis is based on the approximation to the dispersion relation for the Euler equations. This is good for suitable ranges of the parameters measuring dispersive, nonlinear and surface tension effects (and specified below). They also give an idea about how the model may fit real situations. In the context of the approach of \cite{BLS2008} adopted here, the physical regime where the present paper introduces the proposed model is the so-called Benjamin-Ono (BO) regime, theoretically characterized, among others, by the hypothesis of a lower layer of infinite depth but which is useful, as indicated by Kalisch, \cite{Kalisch2007}, to deal with situations where the depth is much larger than the wavelength of a typical wave. 

In addition, the presence of interfacial tension must be under the range of validity specified in \cite{ABR} for the unidirectional model. The asymptotic expansions of the nonlocal operators corresponding to the Benjamin-Ono regime lead to the BO system, introduced in \cite{BLS2008} and investigated in \cite{Xu2012,AnguloS2020,BonaDM2020}, among others. The inclusion of the influence of surface tension at the interface, in the regime specified in \cite{ABR}, leads to the derivation of the two-dimensional asymptotic model proposed in the present paper, that will be called the Benjamin system. This, consequently, becomes the BO system in absence of surface tension. Finally, a similar argument to that considered in \cite{Kalisch2007} to recover the Benjamin equation, after an assumption of uni-directionality of the waves (see \cite{Whitham}), is applied here to the one-dimensional version of the Benjamin system. This leads to a one-parameter family of regularized versions of the Benjamin equation (containing the Benjamin equation as particular case), in a like way to that leading to the regularized BO equation in \cite{KB,BLS2008}.


The rest of the paper is devoted to the analysis of some mathematical properties of the two models introduced, the regularized Benjamin (rBenjamin) equation and the Benjamin system. The study is focused on well-posedness, conserved quantities and existence of solitary wave solutions. First, linear well-posedness is discussed. Furthermore, while the rBenjamin equation is shown to possess at least three functionals preserved by the evolution of smooth enough solutions which decay to zero at infinity as well as a Hamiltonian structure, the Benjamin system only admits linear invariant quantities, but the evolution of candidates for momentum and energy is specified (cf. \cite{BonaDM2020}). Finally, the existence of solitary waves and some of their properties are analyzed in a computational study of comparison of the two models between them and with the Benjamin equation.

The paper is structured as follows. Section \ref{sec:sec2} is devoted to reformulate the Euler system for internal waves with capillary effects in terms of the nonlocal operators used in \cite{BLS2008}. Then the physical regime of validity is incorporated, combining the BO regime with the conditions on the parameter of interfacial tension required by the uni-directional Benjamin equation. The application of these hypotheses leads to the two-dimensional, bi-directional Benjamin system, whose linear well-posedness is studied. From its one-dimensional version, a one-parameter family of the rBenjamin equations is derived. The family contains the usual Benjamin equation as particular case. Linear well-posedness, conserved quantities and Hamiltonian structure of the rBenjamin equation are also analyzed. Existence of solitary wave solutions for the three models and comparisons of the waves are discussed in a computational study in Section \ref{sec:sec3}. Some conclusions and perspectives are outlined in Section \ref{sec:sec4}. 

The following notation will be used throughout the paper. If $s\geq 0, d=1,2$, 
$H^{s}=H^{s}(\mathbb{R}^{d})$ will stand for the $L^{2}-$based Sobolev space of order $s$, with $H^{0}=L^{2}$. The corresponding norm in $H^{s}$ is denoted by $||\cdot||_{s}$. By ${\bf x}$ we will denote the horizontal variable, with ${\bf x}=x$ if $d=1$ and ${\bf x}=(x,y)$ if $d=2$, while $z$ will be used for the vertical variable. The symbol $\nabla_{{\bf x},z}$ (resp. $\Delta_{{\bf x},z}$) will denote the gradient operator (resp. the Laplace operator) with respect to the variables of ${\bf x}$ and $z$ while $\nabla$ (resp. $\Delta$) will only denote the gradient (resp. the Laplacian) with respect to ${\bf x}$. On the other hand, the symbol $\widehat{g}$ will stand for the Fourier transform of $g$ on $\mathbb{R}^{d}$, with ${\bf k}$ as the variable in the Fourier space, with ${\bf k}=k$ if $d=1$ and ${\bf k}=(k_{x},k_{y})$ if $d=2$. Using the notation of \cite{BLS2008}, if $f$ is a function on $\mathbb{R}^{d}$, $f(D)$ will denote the operator defined by the Fourier symbol
\begin{eqnarray*}
\widehat{f(D)g}({\bf k})=f({\bf k})\widehat{g}({\bf k}).
\end{eqnarray*}
In this way, $|D|$ will stand for the operator defined as
\begin{eqnarray*}
\widehat{|D|g}({\bf k})=|{\bf k}|\widehat{g}({\bf k}),
\end{eqnarray*}
with $|{\bf k}|=|k|$ if $d=1$ and $|{\bf k}|=\sqrt{k_{x}^{2}+k_{y}^{2}}$. In the first case, $|D|$ will sometimes be denoted by $\mathcal{H}=\partial_{x}H$, where $H$ is
the Hilbert transform
\begin{eqnarray}
Hf(x)=\frac{1}{\pi}P.V.\int_{-\infty}^{\infty}\frac{f(\xi)}{\xi-x}d\xi,\label{hilbt}
\end{eqnarray}
with $P.V.$ standing for the Principal Value of the integral.

%

%
%
\section{Derivation of the models}
\label{sec:sec2}
\subsection{Euler system for internal waves with surface tension}
\label{sec21}
In this section the full Euler equations for the two-layer interface problem, with nonnegligible surface tension effects at the interface is reformulated by using the approach introduced in \cite{BLS2008}. Let $\zeta=\zeta({\bf x},t)$ denote the deviation of the interface with respect to a rest posed at the vertical variable $z=-d_{1}$ (see Figure \ref{f0}). Assuming that the flows are irrotational, let $\Phi_{i}, i=1,2$ be the corresponding velocity potential for the upper and lower layer respectively. Incompressibility condition means that the potentials satisfy
\begin{eqnarray}
\Delta_{{\bf x},z}\Phi_{i}=0,\; ({\bf x},z)\in\Omega_{t}^{i},\; i=1,2,\label{bens21}
\end{eqnarray}
where the scales are chosen so that the regions $\Omega_{t}^{i}, i=1,2$ occupied by upper and lower layers respectively are described as
\begin{eqnarray*}
\Omega_{t}^{1}&=&\{({\bf x},z)/ -\infty<x,y<\infty, -d_{1}+\zeta({\bf x},t)<z<0\},\\
\Omega_{t}^{2}&=&\{({\bf x},z)/ -\infty<x,y<\infty, -d_{1}-d_{2}<z<-d_{1}+\zeta({\bf x},t)\}.
\end{eqnarray*}
Rigid conditions at the bottom and lid result in the vanishing of the normal component of both velocity potentials at the corresponding boundaries, that is, for $t>0$
\begin{eqnarray}
\Phi_{1z}&=&0,\quad {\rm on}\quad \Gamma_{1}=\{({\bf x},z)/ -\infty<x,y<\infty, z=0\},\label{bens22a}\\
\Phi_{2z}&=&0,\quad {\rm on}\quad \Gamma_{2}=\{({\bf x},z)/ -\infty<x,y<\infty, z=-(d_{1}+d_{2})\}.\label{bens22b}
\end{eqnarray}
On the other hand, the conservation of momentum (Bernoulli equations) for both fluids are
\begin{eqnarray}
\partial_{t}\Phi_{i}+\frac{1}{2}|\nabla_{{\bf x},z}\Phi_{i}|^{2}=-\frac{P_{i}}{\rho_{i}}-gz, \; ({\bf x},z)\in\Omega_{t}^{i},\; i=1,2,\; t>0,\label{bens24}
\end{eqnarray}
where $g$ denotes the acceleration of gravity and $P_{i}$ is the pressure inside the fluid $i, i=1,2$. Finally, the boundary conditions at the interface
\begin{eqnarray*}
 \Gamma_{t}=\{({\bf x},z)/ -\infty<x,y<\infty, z=-d_{1}+\zeta(x,t)\},
\end{eqnarray*}
consist of the assumption that the fluids do not cross the interface (bounding surface condition)
\begin{eqnarray}
\partial_{t}\zeta-\sqrt{1+|\nabla\zeta|^{2}}\partial_{n}\Phi_{i}=0,\quad {\rm on}\quad \Gamma_{t}, \; t\geq 0, \; i=1,2,\label{bens23}
\end{eqnarray}
where $\partial_{n}\Phi_{i}:=\nabla\Phi_{i}\cdot n, i=1,2$, being $n$ the unit upwards normal vector to the interface. Equations (\ref{bens23}) imply that the normal component of the velocity is continuous at the interface, \cite{BLS2008,Lannes}. The boundary conditions are completed with the assumption on continuity of the stress tensor at the interface. In terms of the mean curvature of the deviation
\begin{eqnarray*}
\kappa(\zeta)=-\nabla\cdot\left(\frac{\nabla\zeta}{\sqrt{1+|\nabla\zeta|^{2}}}\right).
\end{eqnarray*}
the continuity condition reads
\begin{eqnarray*}
P_{2}-P_{1}=\sigma\kappa(\zeta),\quad {\rm on}\quad \Gamma_{t},\; t\geq 0, 
\end{eqnarray*}
where $\sigma$ denotes the interfacial tension coefficient.

We now proceed with the reformulation of (\ref{bens21})-(\ref{bens23}) by using the approach given in \cite{BLS2008}. The derivation is a simple adaptation to that given in this paper and therefore we just remind the main points. We introduce the trace of the potentials at the interface
\begin{eqnarray*}
\psi_{i}({\bf x},t)=\Phi_{i}({\bf x},t,-d_{1}+\zeta({\bf x},t)),\; i=1,2.
\end{eqnarray*}
On the other hand, the nonlocal operators considered in the approach are the Dirichlet-to-Neumann (D-N) operator $G[\zeta]$ such that
\begin{eqnarray}
G[\zeta]\psi_{1}=\sqrt{1+|\nabla\zeta|^{2}}\partial_{n}\Phi_{1}\Big|_{z=-d_{1}+\zeta},\label{bens25}
\end{eqnarray}
and the operator $H[\zeta]$ connecting the traces in such a way that
\begin{eqnarray}
H[\zeta]\psi_{1}=\nabla\psi_{2}.\label{bens26}
\end{eqnarray}
The arguments in \cite{BLS2008}, adapted to (\ref{bens21})-(\ref{bens23}), yield the system for $\zeta$ and $\psi_{1}$
\begin{eqnarray}
\partial_{t}\zeta-G[\zeta]\psi_{1}=0,\label{bens27}&&\\
\partial_{t}\left(H[\zeta]\psi_{1}-\nabla\psi_{1}\right)+g(1-\gamma)\nabla\zeta+\frac{1}{2}\nabla\left((H[\zeta]\psi_{1})^{2}-\gamma |\nabla\psi_{1}|^{2}\right)\nonumber &&\\
+\nabla\mathcal{N}(\zeta,\psi_{1})=-\frac{\sigma}{\rho_{2}}\nabla \kappa(\zeta),\label{bens28}&&
\end{eqnarray}
where
\begin{eqnarray}
\mathcal{N}(\zeta,\psi)=\frac{\gamma\left(G[\zeta]\psi+\nabla\zeta\cdot\nabla\psi\right)^{2}-\left(G[\zeta]\psi+\nabla\zeta\cdot H[\zeta]\psi\right)^{2}}{2(\left(\sqrt{1+|\nabla\zeta|^{2}}\right)}.\label{bens29}
\end{eqnarray}

The derivation of the Benjamin system will be made from a dimensionless version of (\ref{bens27})-(\ref{bens29}), where the hypotheses on the physical regime can be applied, and using the same variables as those of \cite{BLS2008}. Thus, if $a$ and $\lambda$ denote, respectively, a typical amplitude and wavelength, we define
\begin{eqnarray*}
\widetilde{{\bf x}}=\frac{{\bf x}}{\lambda},\; \widetilde{z}=\frac{z}{d_{1}},\; \widetilde{t}=\frac{t\sqrt{gd_{1}}}{\lambda},\; \widetilde{\zeta}=\frac{\zeta}{a},\; \widetilde{\psi_{1}}=\frac{1}{a\lambda}\sqrt{\frac{d_{1}}{g}}\psi_{1},
\end{eqnarray*}
and the parameters
\begin{eqnarray}
\epsilon=\frac{a}{d_{1}},\; \mu=\frac{d_{1}^{2}}{\lambda^{2}},\label{bens210}
\end{eqnarray}
with $\gamma=\frac{\rho_{1}}{\rho_{2}}<1, \delta=\frac{d_{1}}{d_{2}}$ denoting the density and depth ratios respectively. The parameters $\epsilon$ and $\mu$ in (\ref{bens210}) represent, respectively, nonlinear and dispersive effects with respect to the upper layer. The corresponding parameters with respect to the lower layer depend on the previous in the form
\begin{eqnarray}
\epsilon_{2}=\frac{a}{d_{2}}=\epsilon\delta,\; \mu_{2}=\frac{d_{2}^{2}}{\lambda^{2}}=\frac{\mu}{\delta^{2}}.\label{bens210b}
\end{eqnarray}
With the arguments used in \cite{BLS2008} the corresponding nondimensional version of (\ref{bens27})-(\ref{bens29}) is given by (tildes are dropped)
\begin{eqnarray}
&&\partial_{t}\zeta-\frac{1}{\mu}G^{\mu}[\epsilon\zeta]\psi_{1}=0,\label{bens211}\\
&&\partial_{t}\left(H^{\mu,\delta}[\epsilon\zeta]\psi_{1}-\nabla\psi_{1}\right)+(1-\gamma)\nabla\zeta
+\frac{\epsilon}{2}\nabla\left((H^{\mu,\delta}[\epsilon\zeta]\psi_{1})^{2}-\gamma |\nabla\psi_{1}|^{2}\right)\nonumber\\
&&+\epsilon\nabla\mathcal{N}^{\mu,\delta}(\epsilon\zeta,\psi_{1})=-\frac{T}{\epsilon\sqrt{\mu}}\nabla \kappa(\epsilon\sqrt{\mu}\zeta),\label{bens212}
\end{eqnarray}
where
\begin{eqnarray*}
\mathcal{N}^{\mu,\delta}(\zeta,\psi)=\mu\frac{\gamma\left(\frac{1}{\mu}G^{\mu}[\zeta]\psi+\nabla\zeta\cdot\nabla\psi\right)^{2}-\left(\frac{1}{\mu}G^{\mu}[\zeta]\psi+\nabla\zeta\cdot H^{\mu,\delta}[\zeta]\psi\right)^{2}}{2(\left(\sqrt{1+\mu|\nabla\zeta|^{2}}\right)},
\end{eqnarray*}
and
\begin{eqnarray}
T=\frac{\sigma}{g\rho_{2}\lambda^{2}}=\frac{\sigma\mu}{g\rho_{2}d_{1}^{2}},\label{bens213}
\end{eqnarray}
(note that $T$ is related to Weber and Froude numbers, cf. \cite{CA,KimA2006})and the nondimensional versions of the operators (\ref{bens25}), (\ref{bens26}) are defined in \cite{BLS2008} as
\begin{eqnarray*}
G^{\mu}[\epsilon\zeta]\psi_{1}&=&-\mu\epsilon\nabla\zeta\cdot\nabla\Phi_{1}\Big|_{z=-1+\epsilon\zeta}+\partial_{z}\Phi_{1}\Big|_{z=-1+\epsilon\zeta},\\
H^{\mu,\delta}[\epsilon\zeta]\psi_{1}&=&\nabla(\Phi_{2}\Big|_{z=-1+\epsilon\zeta}).
\end{eqnarray*}
See \cite{Lannes} for an alternative formulation, based on a system of two equations for $\zeta$ and a combination of the traces $\psi_{i}, i=1,2$.
\subsection{The Benjamin system}
\label{sec22}
In this section we will derive from (\ref{bens211}), (\ref{bens212}) an asymptotic model which is compatible with the physical regime of validation of the Benjamin equation. Thus, we assume that the upper layer is shallow, the deformations are of small amplitude for an infintely deep lower layer and within a Benjamin-Ono regime. In terms of the parameters (\ref{bens210}), (\ref{bens210b}) this means that
\begin{eqnarray}
\mu\sim\epsilon^{2}<<1,\; \mu_{2}=\infty.\label{bens214}
\end{eqnarray}
In addition, and according to the regime associated to the Benjamin model, \cite{ABR}, the parameter (\ref{bens213}) of surface tension at the interface satisfies
\begin{eqnarray}
T\sim\sqrt{\mu}.\label{bens215}
\end{eqnarray}
Conditions (\ref{bens214}) and (\ref{bens215}) determine the asymptotic regime for the Benjamin system. Defining the velocity variable
\begin{eqnarray}
{\bf v}=H^{\mu,\delta}[\epsilon\zeta]\psi_{1}-\gamma\nabla\psi_{1},\label{bens215b}
\end{eqnarray}
we have the asymptotic expansions, \cite{BLS2008}
\begin{eqnarray}
\nabla\psi_{1}&=&-\frac{1}{\gamma}{\bf v}+\frac{\sqrt{\mu}}{\gamma^{2}}|D|{\bf v}+O(\mu),\label{bens216}\\
H^{\mu,\delta}[\epsilon\zeta]\psi_{1}&=&-\sqrt{\mu}|D|\nabla\psi_{1}+O(\mu)=\frac{\sqrt{\mu}}{\gamma}|D|{\bf v}+O(\mu),\label{bens217}\\
\frac{1}{\mu}G^{\mu}[\epsilon\zeta]\psi_{1}&=&\nabla\left((1-\epsilon\zeta)\left(-\frac{1}{\gamma}{\bf v}\right)+\frac{\sqrt{\mu}}{\gamma}|D|{\bf v}\right)+O(\mu),\label{bens218}
\end{eqnarray}
and
\begin{eqnarray}
\frac{T}{\epsilon\sqrt{\mu}}\nabla \kappa(\epsilon\sqrt{\mu}\zeta)=-T\nabla\left(\nabla\cdot\nabla\zeta\right)+O(\epsilon\mu).\label{bens219}
\end{eqnarray}
Using (\ref{bens211}), (\ref{bens218}) and (\ref{bens214}) in the form $\epsilon\sim\sqrt{\mu}$ we have
\begin{eqnarray}
\partial_{t}\zeta=-\frac{1}{\gamma}{\bf v}+O(\epsilon).\label{bens220}
\end{eqnarray}
Now, introducing the modelling parameter $\alpha\geq 0$ as in \cite{BLS2008} and using (\ref{bens220}) lead to
\begin{eqnarray}
\nabla\cdot{\bf v}=(1-\alpha)\nabla\cdot{\bf v}+\alpha\nabla\cdot{\bf v}=(1-\alpha)\nabla\cdot{\bf v}-\alpha\gamma \partial_{t}\zeta+O(\epsilon).\label{bens221}
\end{eqnarray}
Finally, applying (\ref{bens216})-(\ref{bens221}) to (\ref{bens211}), (\ref{bens212}) and dropping the $O(\epsilon^{2})$ terms we obtain the following system for the deviation of the interface $\zeta$ and the variable ${\bf v}$ in (\ref{bens215b})
\begin{eqnarray}
\left(1+\frac{\alpha\sqrt{\mu}}{\gamma}|D|\right)\partial_{t}\zeta+\frac{1}{\gamma}\nabla\left((1-\epsilon\zeta){\bf v}\right)-(1-\alpha)\frac{\sqrt{\mu}}{\gamma^{2}}|D|\nabla\cdot{\bf v}=0,&&\label{bens222}\\
\partial_{t}{\bf v}+(1-\gamma)\nabla\zeta-\frac{\epsilon}{2\gamma}\nabla\left(|{\bf v}|^{2}\right)=T\nabla(\nabla\cdot\nabla\zeta),&&\label{bens223}
\end{eqnarray}
where $T$ is given by (\ref{bens213}).
\begin{remark}
When interfacial tension is negligible, then $T=0$ and (\ref{bens222}), (\ref{bens223}) is the BO system derived in \cite{BLS2008}.  Furthermore, this and (\ref{bens219}) imply that the internal wave equations (\ref{bens211}), (\ref{bens212}) are consistent with the Benjamin system (\ref{bens222}), (\ref{bens223}), in the sense defined in \cite{BLS2008}, with a precision $O(\mu)$ (see Theorem 6 in that reference).
\end{remark}
We analyze now some elementary mathematical properties of (\ref{bens222}), (\ref{bens223}). Concerning linear well-posedness, we consider the associated linear problem
\begin{eqnarray}
\left(1+\frac{\alpha\sqrt{\mu}}{\gamma}|D|\right)\partial_{t}\zeta+\frac{1}{\gamma}\nabla\left({\bf v}\right)-(1-\alpha)\frac{\sqrt{\mu}}{\gamma^{2}}|D|\nabla\cdot{\bf v}&=&0,\label{bens32a}\\
\partial_{t}{\bf v}+(1-\gamma)\nabla\zeta&=&T\nabla(\nabla\cdot\nabla\zeta).\label{bens32b}
\end{eqnarray}
Note that the operator 
\begin{eqnarray}
J(\alpha)=1+\frac{\alpha\sqrt{\mu}}{\gamma}|D|,\label{joper}
\end{eqnarray} 
has Fourier symbol
\begin{eqnarray*}
\widehat{J(\alpha)}({\bf k})=1+\frac{\alpha\sqrt{\mu}}{\gamma}|{\bf k}|.
\end{eqnarray*}
 It is therefore invertible for $\alpha\geq 0$ and we can write (\ref{bens32a}), (\ref{bens32b}) in the form
\begin{eqnarray}
\partial_{t}\zeta+\frac{1}{\gamma}J(\alpha)^{-1}J(\alpha-1)\nabla\left({\bf v}\right)&=&0,\label{bens33a}\\
\partial_{t}{\bf v}+\nabla\left((1-\gamma)-T\nabla\cdot\nabla\right)\zeta&=&0.\label{bens33b}
\end{eqnarray}
We now take Fourier transform in (\ref{bens33a}), (\ref{bens33b}) with respect to the spatial variables to have, for ${\bf v}=(v_{1},v_{2})$
\begin{eqnarray}
\frac{d}{dt}\begin{pmatrix}\widehat{\zeta}({\bf k},t)\\ \widehat{v_{1}}({\bf k},t)\\
\widehat{v_{2}}({\bf k},t)\end{pmatrix}+i|{\bf k}|\mathcal{A}({\bf k})\begin{pmatrix}\widehat{\zeta}({\bf k},t)\\ \widehat{v_{1}}({\bf k},t)\\
\widehat{v_{2}}({\bf k},t)\end{pmatrix}=0,\label{bens34}
\end{eqnarray}
where
\begin{eqnarray*}
\mathcal{A}({\bf k})=\begin{pmatrix}0&\frac{k_{x}}{\gamma|{\bf k}|}\left(\frac{\widehat{J(\alpha-1)}({\bf k})}{\widehat{J(\alpha)}({\bf k})}\right)&\frac{k_{y}}{\gamma|{\bf k}|}\left(\frac{\widehat{J(\alpha-1)}({\bf k})}{\widehat{J(\alpha)}({\bf k})}\right)\\
\frac{k_{x}}{\gamma|{\bf k}|}\left(1-\gamma+T|{\bf k}|^{2}\right)&0&0\\
\frac{k_{y}}{\gamma|{\bf k}|}\left(1-\gamma+T|{\bf k}|^{2}\right)&0&0\end{pmatrix}
\end{eqnarray*}
The eigenvalues of $\mathcal{A}({\bf k})$ are $\{0,\pm\sigma({\bf k})\}$ with
\begin{eqnarray*}
\sigma({\bf k})=\left(\frac{\left(1-\gamma+T|{\bf k}|^{2}\right)}{\gamma}\left(\frac{\widehat{J(\alpha-1)}({\bf k})}{\widehat{J(\alpha)}({\bf k})}\right)\right)^{1/2},
\end{eqnarray*}
which implies linear well-posedness of  (\ref{bens222}), (\ref{bens223}) whenever $\alpha\geq 1$. In order to specify the Sobolev spaces in detail, we diagonalize the matrix $\mathcal{A}({\bf k})$ (see e.~g. \cite{DougalisMS2007})
\begin{eqnarray*}
P({\bf k})^{-1}\mathcal{A}({\bf k})P({\bf k})=\begin{pmatrix} 0&0&0\\0&\sigma({\bf k})&0\\0&0&-\sigma({\bf k})\end{pmatrix}
\end{eqnarray*}
with
\begin{eqnarray*}
P({\bf k})=\begin{pmatrix} 0&\tau({\bf k})&-\tau({\bf k})\\-\frac{k_{y}}{\gamma|{\bf k}|}&\frac{k_{x}}{\gamma|{\bf k}|}&\frac{k_{x}}{\gamma|{\bf k}|}\\
\frac{k_{x}}{\gamma|{\bf k}|}&\frac{k_{y}}{\gamma|{\bf k}|}&\frac{k_{x}}{\gamma|{\bf k}|}\end{pmatrix},\;
P({\bf k})^{-1}=\frac{1}{2\tau({\bf k})}\begin{pmatrix} 0&-2\tau({\bf k})\frac{k_{y}}{\gamma|{\bf k}|}&2\tau({\bf k})\frac{k_{x}}{\gamma|{\bf k}|}\\1&\tau({\bf k})\frac{k_{x}}{\gamma|{\bf k}|}&\tau({\bf k})\frac{k_{y}}{\gamma|{\bf k}|}\\
-1&\tau({\bf k})\frac{k_{x}}{\gamma|{\bf k}|}&\tau({\bf k})\frac{k_{y}}{\gamma|{\bf k}|}\end{pmatrix},
\end{eqnarray*}
and
\begin{eqnarray*}
\tau({\bf k})=\frac{\sigma({\bf k})}{\left(1-\gamma+T|{\bf k}|^{2}\right)}.
\end{eqnarray*}
With the change of variables
\begin{eqnarray}
\begin{pmatrix} \widehat{\eta}\\\widehat{{\bf w}}\end{pmatrix}=P^{-1}\begin{pmatrix} \widehat{\zeta}\\\widehat{{\bf v}}\end{pmatrix},\; {\bf w}=(w_{1},w_{2}),\label{bens35}
\end{eqnarray}
the system (\ref{bens34}) is transformed into 
\begin{eqnarray*}
\frac{d}{dt}\begin{pmatrix}\widehat{\eta}({\bf k},t)\\ \widehat{w_{1}}({\bf k},t)\\
\widehat{w_{2}}({\bf k},t)\end{pmatrix}+i|{\bf k}|\begin{pmatrix} 0&0&0\\0&\sigma({\bf k})&0\\0&0&-\sigma({\bf k})\end{pmatrix}\begin{pmatrix}\widehat{\eta}({\bf k},t)\\ \widehat{w_{1}}({\bf k},t)\\
\widehat{w_{2}}({\bf k},t)\end{pmatrix}=0,
\end{eqnarray*}
with solution
\begin{eqnarray*}
\widehat{\eta}({\bf k},t)=\widehat{\eta}({\bf k},0),\; 
\widehat{w_{1}}({\bf k},t)=e^{-i|{\bf k}|\sigma({\bf k})t}\widehat{w_{1}}({\bf k},0),\;
\widehat{w_{2}}({\bf k},t)=e^{i|{\bf k}|\sigma({\bf k})t}\widehat{w_{2}}({\bf k},0).
\end{eqnarray*}
Now, since
\begin{eqnarray}
\widehat{\eta}({\bf k})&=&-\frac{k_{y}}{|{\bf k}|}\widehat{v_{1}}({\bf k})+\frac{k_{x}}{|{\bf k}|}\widehat{v_{2}}({\bf k}),\nonumber\\
\widehat{w_{1}}({\bf k})&=&\frac{1}{2\tau({\bf k})}\widehat{\zeta}({\bf k})+\frac{k_{x}}{2|{\bf k}|}\widehat{v_{1}}({\bf k})+\frac{k_{y}}{2|{\bf k}|}\widehat{v_{2}}({\bf k}),\nonumber\\
\widehat{w_{2}}({\bf k})&=&-\frac{1}{2\tau({\bf k})}\widehat{\zeta}({\bf k})+\frac{k_{x}}{2|{\bf k}|}\widehat{v_{1}}({\bf k})+\frac{k_{y}}{2|{\bf k}|}\widehat{v_{2}}({\bf k}),\label{bens35b}
\end{eqnarray}
and $\tau({\bf k})$ has order $-1$, then
\begin{eqnarray*}
(\zeta,v_{1},v_{2})\in H^{s+1}\times H^{s}\times H^{s}\Rightarrow (\eta,w_{1},w_{2})\in H^{s}\times H^{s}\times H^{s},\; s> 0,
\end{eqnarray*}
and therefore, when $\alpha\geq 1$ the system (\ref{bens222}), (\ref{bens223}) is linearly well-posed in $H^{s+1}\times H^{s}\times H^{s}, s>0$.

Another point of interest is the existence of conserved quantities for the one-dimensional version (in the $x-$direction) of  (\ref{bens222}), (\ref{bens223}) 
\begin{eqnarray}
\left(1+\frac{\alpha\sqrt{\mu}}{\gamma}\mathcal{H}\right)\partial_{t}\zeta+\frac{1}{\gamma}\partial_{x}\left((1-\epsilon\zeta){u}\right)-(1-\alpha)\frac{\sqrt{\mu}}{\gamma^{2}}\mathcal{H}\partial_{x}{u}&=&0,\label{bens224}\\
\partial_{t}{u}+(1-\gamma)\nabla\zeta-\frac{\epsilon}{2\gamma}\nabla\left(u^{2}\right)&=&T\partial_{x}^{3}\zeta,\label{bens225}
\end{eqnarray}
where $\mathcal{H}=\partial_{x}H$, being $H$ the Hilbert transform (\ref{hilbt}),
and $u$ is a horizontal velocity-like variable. The one-dimensional version (\ref{bens224}), (\ref{bens225})  trivially admits the linear functionals
\begin{eqnarray*}
I_{1}(\zeta,u)=\int_{-\infty}^{\infty}\zeta dx,\quad I_{2}(\zeta,u)=\int_{-\infty}^{\infty}u dx,
\end{eqnarray*}
as invariants by the time evolution of smooth enough solutions which decay, along with higher-order derivatives, to zero at infinity. As in the case of the BO system, \cite{BonaDM2020}, no other conserved quantities were found. It may be worth mentioning that the quantities
\begin{eqnarray*}
I&=&\int_{-\infty}^{\infty}\zeta u dx,\\
H&=&\frac{1-\gamma}{2}\int_{-\infty}^{\infty}\zeta^{2} dx-\frac{\varepsilon}{2\gamma}\int_{-\infty}^{\infty}\zeta u^{2} dx\\
&&+\frac{1}{2\gamma}\int_{-\infty}^{\infty}u^{2} dx
-(1-\alpha)\frac{\sqrt{\mu}}{2\gamma^{2}}\int_{-\infty}^{\infty}u\mathcal{H}u dx-\frac{T}{2}\int_{-\infty}^{\infty}\zeta_{x}^{2}dx,
\end{eqnarray*}
satisfy the time evolution given by the equations (cf. \cite{BonaDM2020})
\begin{eqnarray*}
\frac{d}{dt}I=-\frac{\alpha\sqrt{\mu}}{\gamma}\int_{-\infty}^{\infty}u\mathcal{H}\zeta_{t} dx 
\end{eqnarray*}
\begin{eqnarray*}
\frac{d}{dt}H=-\frac{\alpha\sqrt{\mu}}{\gamma}\int_{-\infty}^{\infty}\left((1-\gamma)\zeta-\frac{\varepsilon}{2\gamma}u^{2}\right)\mathcal{H}\zeta_{t} dx.
\end{eqnarray*}
Thus they are preserved only in the case $\alpha=0$.
\subsection{The regularized Benjamin equation}
\label{sec23}
The reduction of a two-way model like (\ref{bens224}), (\ref{bens225}) to corresponding unidirectional models can be formally made in the same way as in \cite{BonaDM2020} for the Benjamin-Ono and Intermediate Long Wave systems for internal waves or as that of the bidirectional model for interfacial capillary-gravity waves in deep water in \cite{Kalisch2007}. Summarizing (see \cite{Whitham}), when $O(\epsilon)=O(\sqrt{\mu})$ terms in (\ref{bens224}), (\ref{bens225}) are neglected, the deformation $\zeta$ will satisfy a wave equation with speed
\begin{eqnarray}
c_{\gamma}=\sqrt{\frac{1-\gamma}{\gamma}}.\label{bens226}
\end{eqnarray} 
Right moving waves will then be of the form $\zeta=\zeta_{0}(x-c_{\gamma}t), u=\sqrt{\gamma(1-\gamma)}\zeta$. In order to find solutions of (\ref{bens224}), (\ref{bens225}) moving e.~g. to the right to order $O(\epsilon)=O(\sqrt{\mu})$ we assume $u$ of the form
\begin{eqnarray}
u=\sqrt{\gamma(1-\gamma)}\left(\zeta+A\epsilon+B\sqrt{\mu}\right),\label{bens227}
\end{eqnarray}
for some functions $A, B$ of $\zeta$. As in \cite{BonaDM2020}, we substitute (\ref{bens227}) into (\ref{bens224}), (\ref{bens225}) and retain only the $O(\epsilon)=O(\sqrt{\mu})$ terms. Then consistency of the two equations leads to
\begin{eqnarray}
A=\frac{1}{4}\zeta^{2},\; B=\frac{1}{2\gamma}\mathcal{H}\zeta-\frac{T}{2(1-\gamma)\sqrt{\mu}}\partial_{x}^{2}\zeta,\label{bens228}
\end{eqnarray} 
and the one-parameter family of unidirectional equations for $\zeta$
\begin{eqnarray}
\left(1+\frac{\alpha\sqrt{\mu}}{\gamma}\mathcal{H}\right)\partial_{t}\zeta+c_{\gamma}\partial_{x}\zeta-\frac{3\epsilon}{4}c_{\gamma}\partial_{x}\zeta^{2}-c_{\gamma}\frac{(1-2\alpha)}{2\gamma}\sqrt{\mu}\mathcal{H}\partial_{x}\zeta&&\nonumber\\
-\frac{T}{2\sqrt{\gamma(1-\gamma)}}\partial_{x}^{3}\zeta=0.&&\label{bens229}
\end{eqnarray}
In absence of surface tension ($T=0$) equation (\ref{bens229}) reduces to the regularized Benjamin-Ono equation, \cite{KB,BLS2008}. According to this, (\ref{bens229}) will be called the regularized Benjamin (or rBenjamin) equation. The Benjamin equation corresponds to $\alpha=0$.

Taking the Fourier transform in the corresponding linearized equation of (\ref{bens229}) leads to solutions of the form $\widehat{\zeta}(k,t)=e^{-ikm(k)t}\widehat{\zeta}(k,0)$ with
\begin{eqnarray*}
m(k)=\frac{c_{\gamma}\left(1-\frac{(1-2\alpha)}{2\gamma}\sqrt{\mu}|k|+\widetilde{T}k^{2}\right)}{1+\frac{\alpha}{\gamma}\sqrt{\mu}|k|},\; \widetilde{T}=\frac{T}{2\sqrt{\gamma(1-\gamma)}},
\end{eqnarray*}
and to the linear dispersion relation $
\omega(k)=km(k)$, ensuring linear well-posedness of (\ref{bens229}).

In analogy with the regularized BO equation, \cite{BonaDM2020}, (\ref{bens229}) admits at least three time invariant functionals and a Hamiltonian structure. The conserved quantities are
\begin{eqnarray*}
&& C(\zeta)=\int_{-\infty}^{\infty} \zeta dx,\qquad  D(\zeta)=\frac{1}{2}\int_{-\infty}^{\infty} \left(\zeta^{2}+\sqrt{\mu}\frac{\alpha}{\gamma}\zeta\mathcal{H}\zeta\right) dx,\\
&& E(\zeta)=\frac{c_{\gamma}}{2}\int_{-\infty}^{\infty} \left(\zeta^{2}-\sqrt{\mu}\frac{(1-2\alpha)}{2\gamma}\zeta\mathcal{H}\zeta-\frac{1}{2}\zeta^{3}\right) dx+\frac{\widetilde{T}}{2}\int_{-\infty}^{\infty}(\partial_{x}\zeta)^{2}dx.
\end{eqnarray*}
The last one enables (\ref{bens229}) to have a Hamiltonian formulation
\begin{eqnarray*}
\partial_{t}\zeta=\mathcal{J}\frac{\delta}{\delta\zeta}E(\zeta),
\end{eqnarray*}
with structure operator $\mathcal{J}=-\partial_{x}J(\alpha)$, where $J(\alpha)$ is given by (\ref{joper}) (in its one-dimensional version) and $\frac{\delta}{\delta\zeta}$ denotes variational (Fr\'echet) derivative. When $T=0$ we recover the invariants and Hamiltonian structure of the regularized BO equation, \cite{BonaDM2020}. As in the particular case of the Benjamin equation ($\alpha=0$), the existence of these conserved quantities and the theory of Kenig et al., \cite{KenigPV1991,KenigPV1993}, might be used to obtain local and global well-posedness results for (\ref{bens229}), see \cite{Linares1999,LinaresS2005}.


\section{A computational study of solitary wave solutions}
\label{sec:sec3}
\subsection{Preliminaries}
\label{sec31}
Another property typically studied in water wave models is the existence of special solutions. Of particular interest are the solutions of solitary-wave type, due to their relevance in the general dynamics of some models, \cite{Bona1981}. The present section is concerned with this topic for the cases of the one-dimensional version of the Benjamin system (\ref{bens224}), (\ref{bens225}), and the regularized Benjamin equations (\ref{bens229}). The purpose here, developed by computational means, is two-fold: The first one is related to the existence of solitary-wave solutions, for which no theoretical results are available. On the other hand, we are also interested in comparing these solitary wave profiles of the two models and with those of the Benjamin equation. In particular, and following here the study developed for the system introduced in \cite{Kalisch2007}, the dynamics of solitary wave solutions of the Benjamin equation under the evolution given by both the Benjamin system and the rBenjamin equation is numerically investigated.

We start with a description of the equations involved for the generation of solitary wave solutions and the numerical tools used for the computational study. For the case of the Benjamin system (\ref{bens224}), (\ref{bens225}), we are looking for solutions in the form of traveling waves $\zeta(x,t)=\zeta(x-c_{s}t), u(x,t)=u(x-c_{s}t)$, for some speed $c_{s}\neq 0$ and with the profiles $\zeta(X), u(X), X=x-c_{s}t$ which are smooth and decay to zero as $|X|\rightarrow\infty$. Substituting into (\ref{bens224}), (\ref{bens225}) and integrating once yield the system
\begin{eqnarray}
\begin{pmatrix}-c_{s}J(\alpha)&\frac{1}{\gamma}J(\alpha-1)\\(1-\gamma)-T\partial_{x}^{2}&-c_{s}\end{pmatrix}\begin{pmatrix}\zeta\\u\end{pmatrix}=\frac{\epsilon}{\gamma}\begin{pmatrix}\zeta u\\ \frac{u^{2}}{2}\end{pmatrix}.\label{bens41}
\end{eqnarray}
In the case of the rBenjamin equation (\ref{bens229}), solutions $\zeta(x,t)=\zeta(x-c_{s}t), c_{s}>0$, with smooth profiles $\zeta(X)\rightarrow 0, |X|\rightarrow\infty$ will satisfy the equation
\begin{eqnarray}
-c_{s}J(\alpha)\zeta+c_{\gamma}\left(\zeta-\frac{3\epsilon}{4}\zeta^{2}+\frac{2\alpha-1}{2\gamma}\sqrt{\mu}\mathcal{H}\zeta\right)-T\zeta^{\prime\prime}=0.\label{bens42}
\end{eqnarray}
We will focus here on the strong numerical evidence of existence of solutions of (\ref{bens41}) and (\ref{bens42}) and on the properties of the solitary waves suggested by the computations, leaving the study of theoretical results to some future research. The computational approach is based on the iterative resolution of the algebraic systems obtained from the Fourier representation of (\ref{bens41}) and (\ref{bens42}). Written in fixed-point form, these are respectively
\begin{eqnarray}
\begin{pmatrix}-c_{s}\widehat{J(\alpha)}(k)&\frac{1}{\gamma}\widehat{J(\alpha-1)}(k)\\(1-\gamma)+Tk^{2}&-c_{s}\end{pmatrix}\begin{pmatrix}\widehat{\zeta}(k)\\\widehat{u}(k)\end{pmatrix}=\frac{\epsilon}{\gamma}\begin{pmatrix}\widehat{\zeta u}(k)\\ \widehat{\frac{u^{2}}{2}}(k)\end{pmatrix}.\label{bens43}
\end{eqnarray}
for $k\in\mathbb{R}$, and where $\widehat{J(r)}(k)=1+\frac{r}{\gamma}\sqrt{\mu}|k|$, and
\begin{eqnarray}
\left(-c_{s}\widehat{J(\alpha)}(k)+c_{\gamma}\left(1+\frac{2\alpha-1}{2\gamma}\sqrt{\mu}|k|\right)+Tk^{2}\right)\widehat{\zeta}(k)=\frac{3c_{\gamma}\epsilon}{4}\widehat{\zeta^{2}}(k).\label{bens44}
\end{eqnarray}
The numerical procedure is performed in the standard way: for each case (\ref{bens41}), (\ref{bens42}), the corresponding periodic problem on a long enough interval $(-L,L)$ is implemented via the Fourier representation based on the form (\ref{bens43}) and (\ref{bens44}) respectively, where now $k\in\mathbb{Z}$ and $\widehat{\zeta}(k), \widehat{u}(k)$ represent the corresponding $k$th Fourier coefficient. The resulting algebraic equations for each $k$ are iteratively solved by using the Petviashvili's method, \cite{Petv1976,pelinovskys}, taking advantage of the homogeneous character (of degree two) of the nonlinear term. The Petviashvili's iteration is complemented with vector extrapolation techniques, \cite{sidi,sidifs,smithfs}, in order to accelerate the convergence, \cite{AlvarezD2015}.
\subsection{A comparative study}
The approximate solitary wave solutions of the Benjamin system (\ref{bens224}), (\ref{bens225}), the regularized Benjamin equation (\ref{bens229}) and the Benjamin equation (again (\ref{bens229}) with $\alpha=0$) are here compared in a series of numerical experiments. The dimensionless parameters for the computations are taken as $\epsilon=\sqrt{\mu}=0.1$, with $T=0.1$, while different values of $\gamma$ and $c_{s}$ are considered. The interval of approximation is determined by $L=256$ and $4096$ Fourier modes are typically used.
\begin{figure}[!htbp] 
\centering
{\includegraphics[width=\columnwidth]{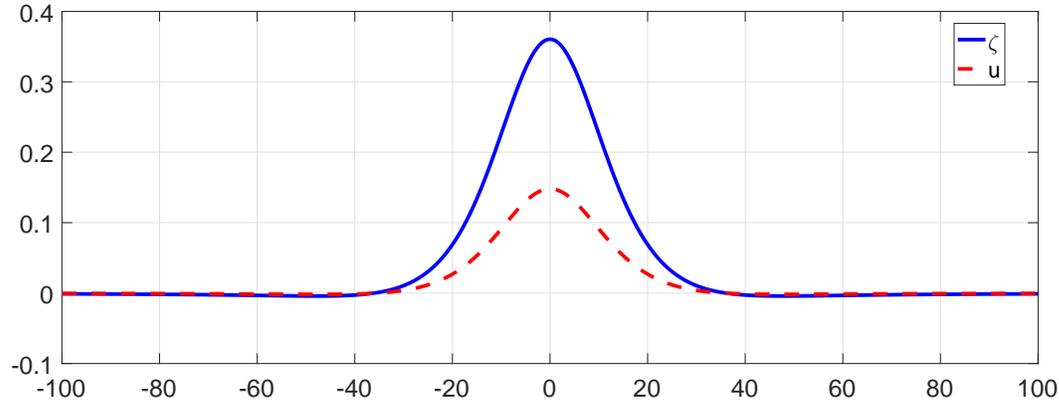}}
\caption{Approximate profiles $(\zeta,u)$ of Benjamin system (\ref{bens41}) with $\gamma=0.8, c_s=0.49$.
}
\label{fig:bens_fig1}
\end{figure}
\begin{figure}[!htbp] 
\centering
{\includegraphics[width=\columnwidth]{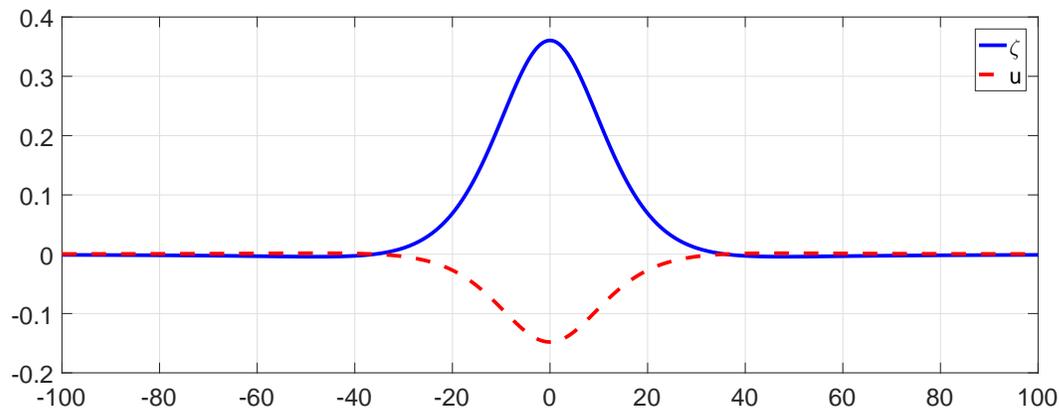}}
\caption{Approximate profiles $(\zeta,u)$ of Benjamin system (\ref{bens41}) with $\gamma=0.8, c_s=-0.49$.
}
\label{fig:bens_fig2}
\end{figure}

Note first that if $(c_{s},\zeta,u)$ is a solution of (\ref{bens41}) then $(-c_{s},\zeta,-u)$ is also a solution, with the same $\zeta$ profile traveling in opposite direction. This is illustrated in Figures \ref{fig:bens_fig1} and \ref{fig:bens_fig2}, with the representation of the approximate $\zeta$ and $u$ solitary wave profiles corresponding to $\gamma=0.8$ and $c_{s}=0.49$.

\begin{figure}[htbp]
\centering
\subfigure[]
{\includegraphics[width=\columnwidth]{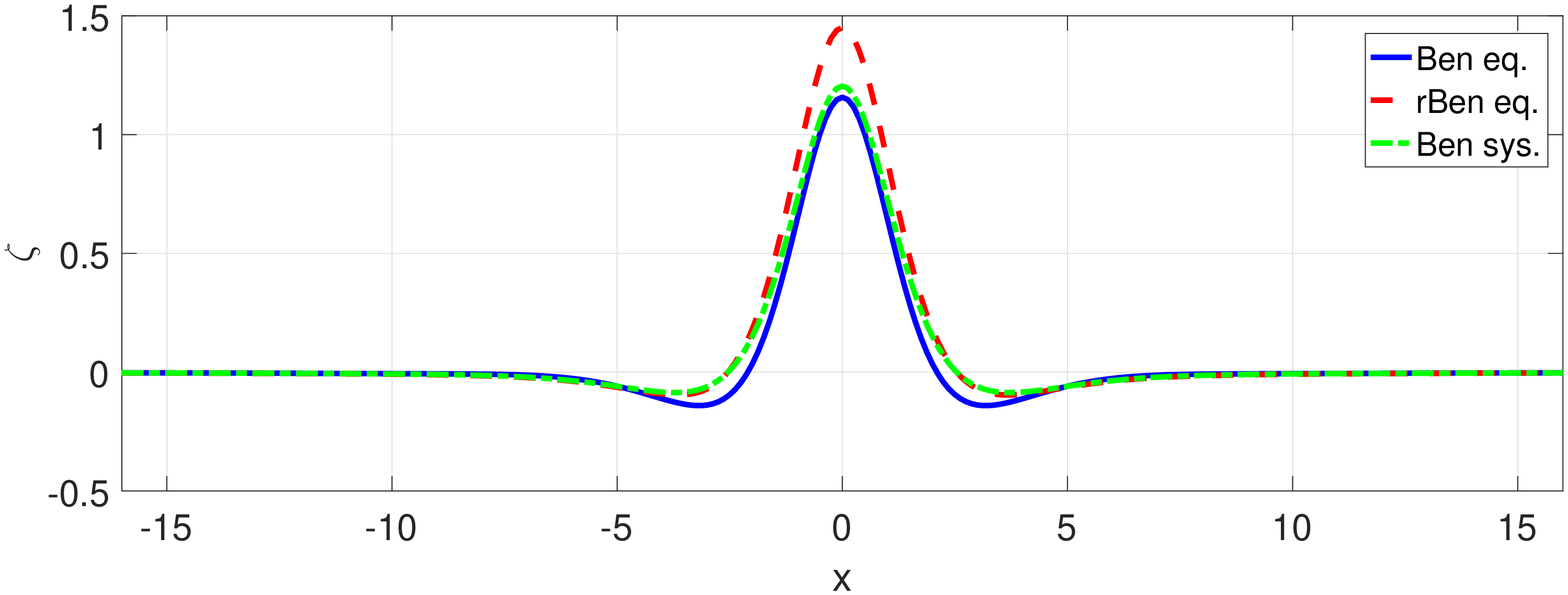}}
\subfigure[]
{\includegraphics[width=\columnwidth]{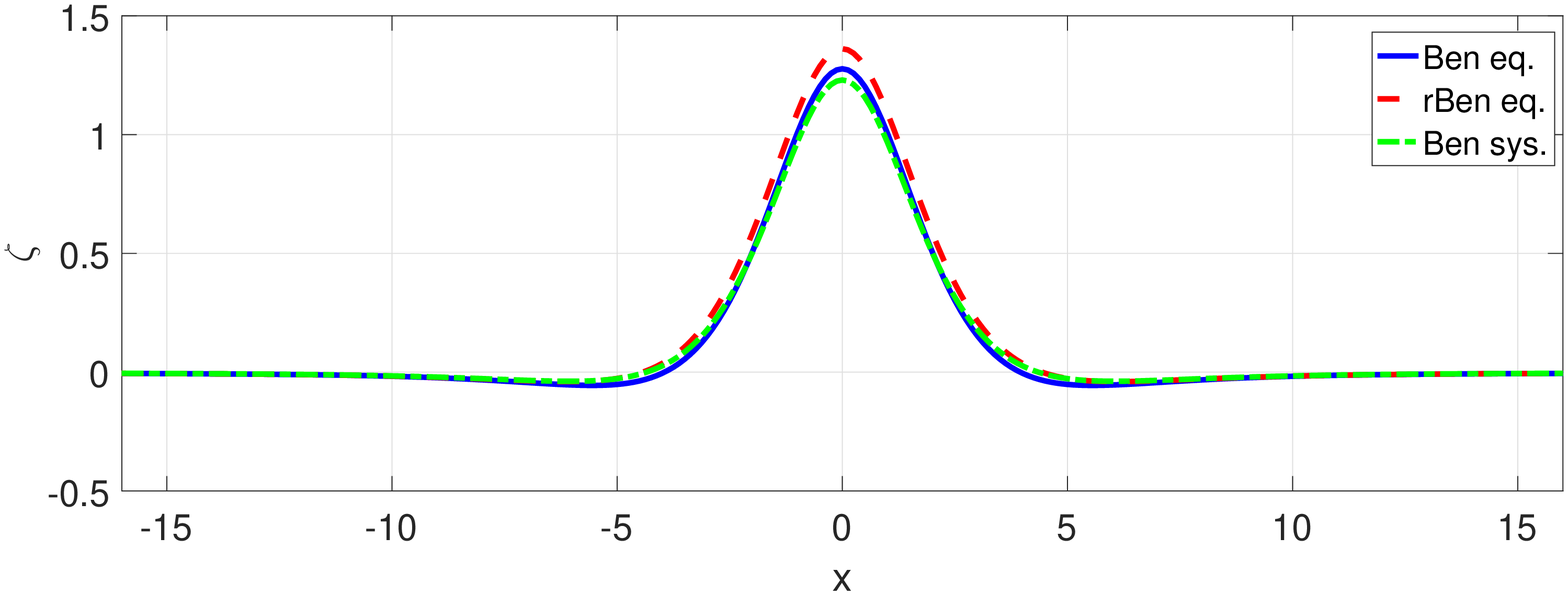}}
\subfigure[]
{\includegraphics[width=\columnwidth]{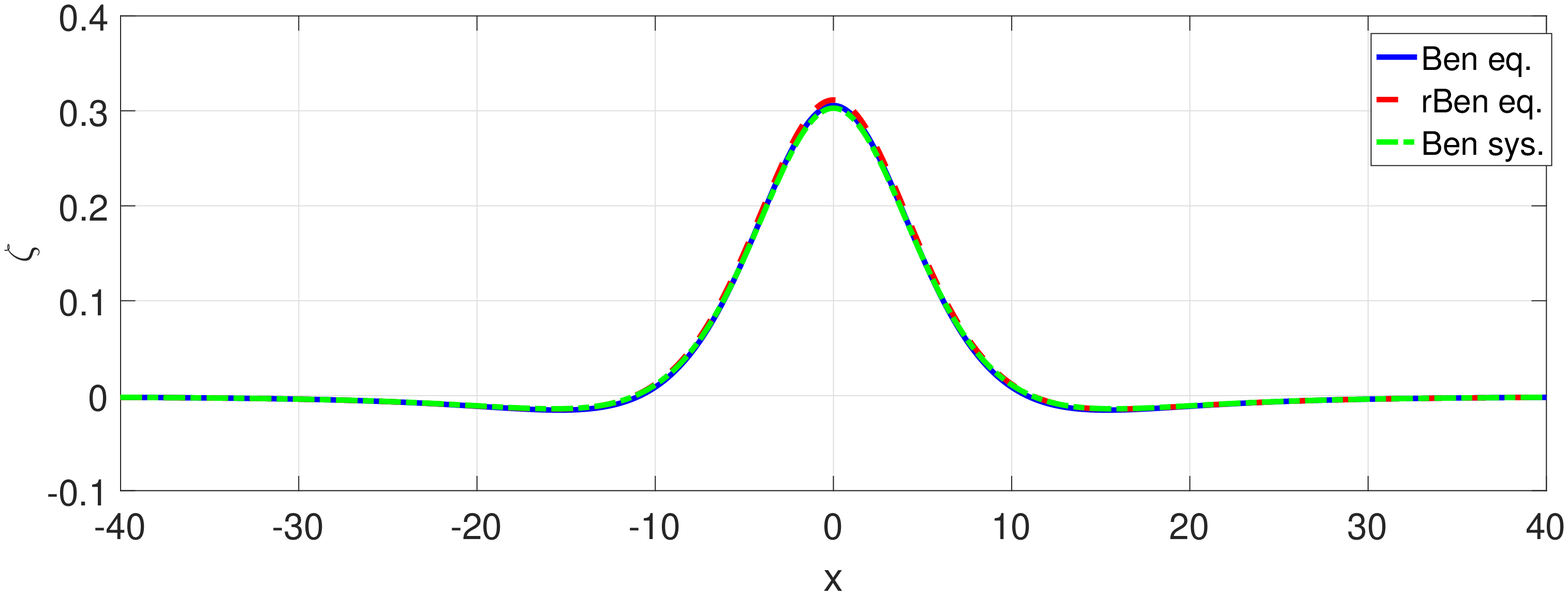}}
\caption{Comparison of approximate $\zeta$ solitary wave profiles. (a) $\gamma=0.4, c_{s}=1.1$; (b) $\gamma=0.6, c_{s}=0.75$; (c) $\gamma=0.8, c_{s}=0.49$.}
\label{fig:bens_fig3}
\end{figure}
\begin{figure}[htbp]
\centering
\subfigure[]
{\includegraphics[width=\columnwidth]{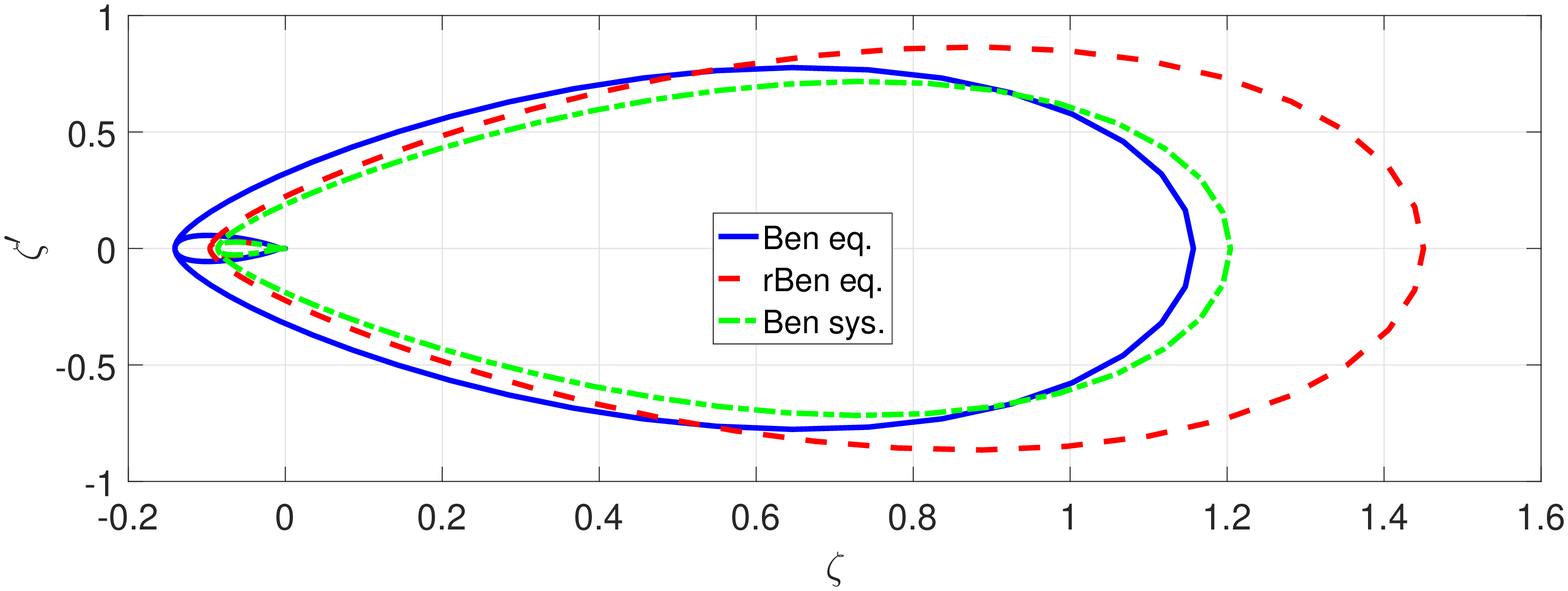}}
\subfigure[]
{\includegraphics[width=\columnwidth]{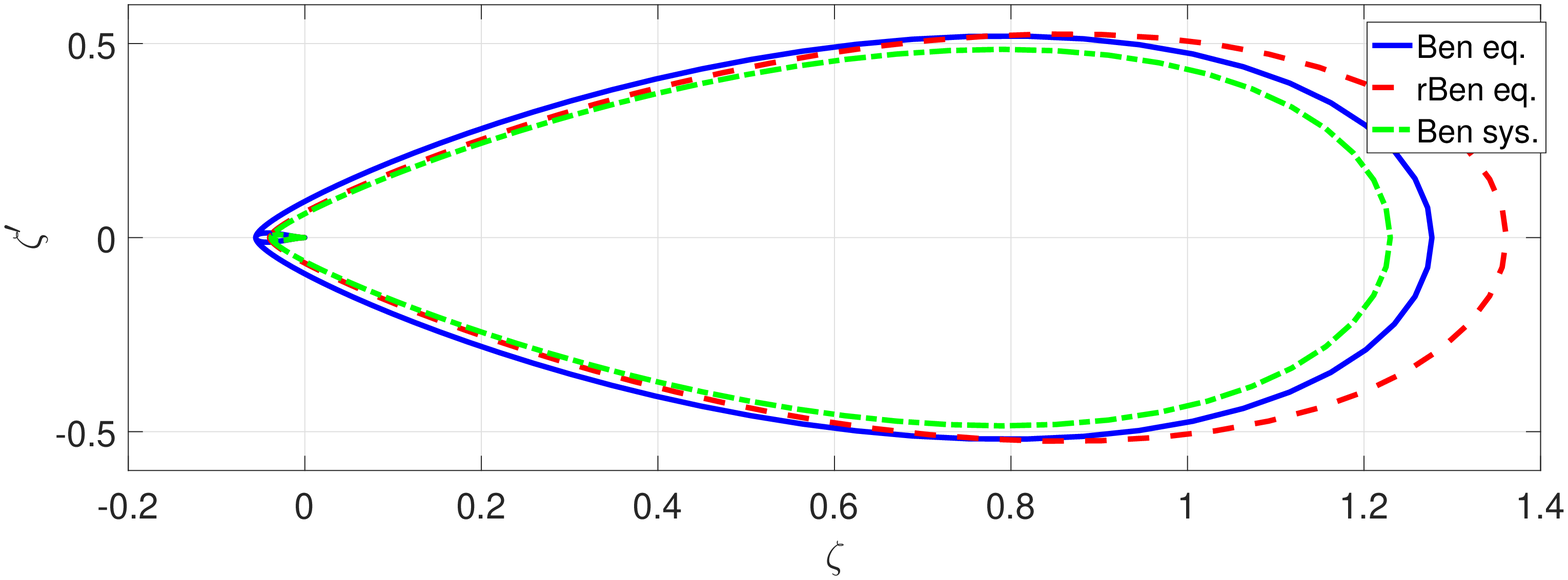}}
\subfigure[]
{\includegraphics[width=\columnwidth]{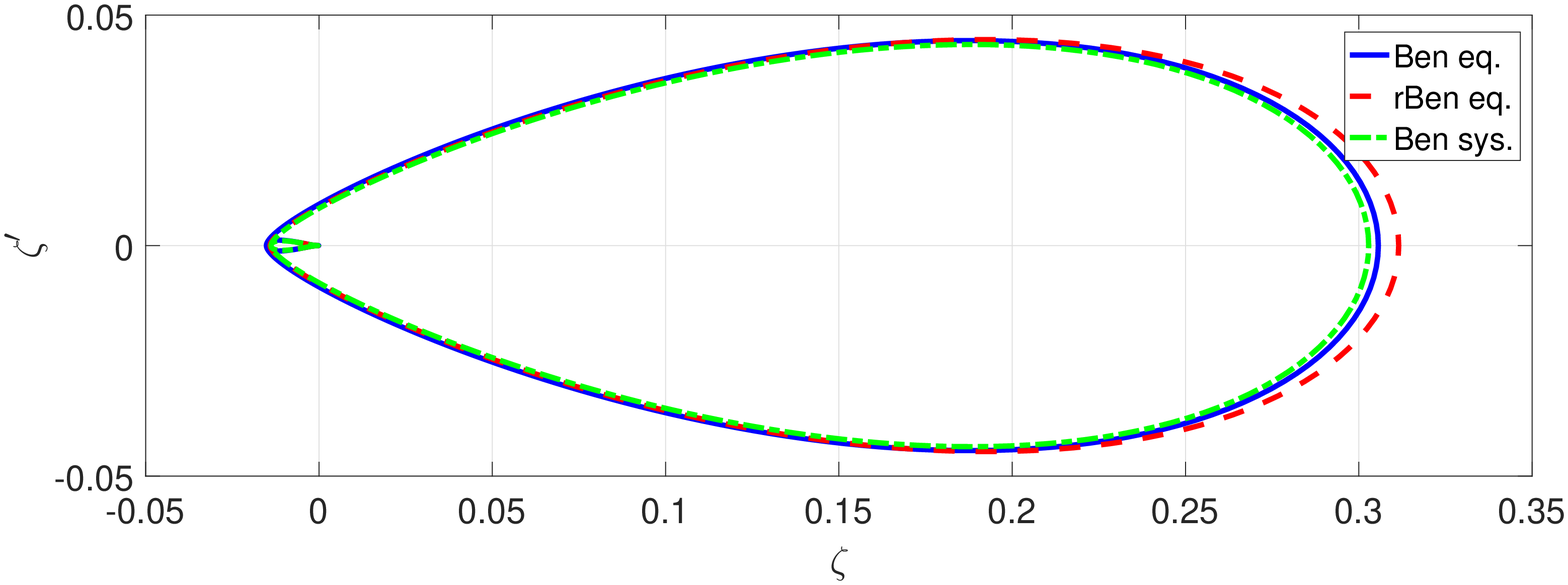}}
\caption{Phase portraits of the approximate $\zeta$ profiles of Figure \ref{fig:bens_fig3}.}
\label{fig:bens_fig4}
\end{figure}
The second observation is that the computations generate solitary wave profiles when $|c_{s}|<c_{\gamma}$, where $c_{\gamma}$ is given by (\ref{bens226}), and the $\zeta$ profiles are of elevation. (This cannot serve us, indeed, to discard the existence of solitary wave solutions of depression as in the BO system, \cite{BonaDM2020,AnguloS2020}.) Additional observations are that the profiles are not positive; they contain an oscillatory decay in the same way as the known behaviour of the solitary wave solutions of the Benjamin equation, \cite{ABR,Benjamin1992}. The waves are taller and with less oscillations as $|c_{s}|$ moves away from $c_{\gamma}$. These properties may be compared with those of the solitary wave solutions of the BO and ILW systems, studied theoretically in \cite{AnguloS2020} and computationally in \cite{BonaDM2020}, for which the solitary waves do not show oscillatory decay.

\begin{figure}[htbp]
\centering
\subfigure[]
{\includegraphics[width=\columnwidth]{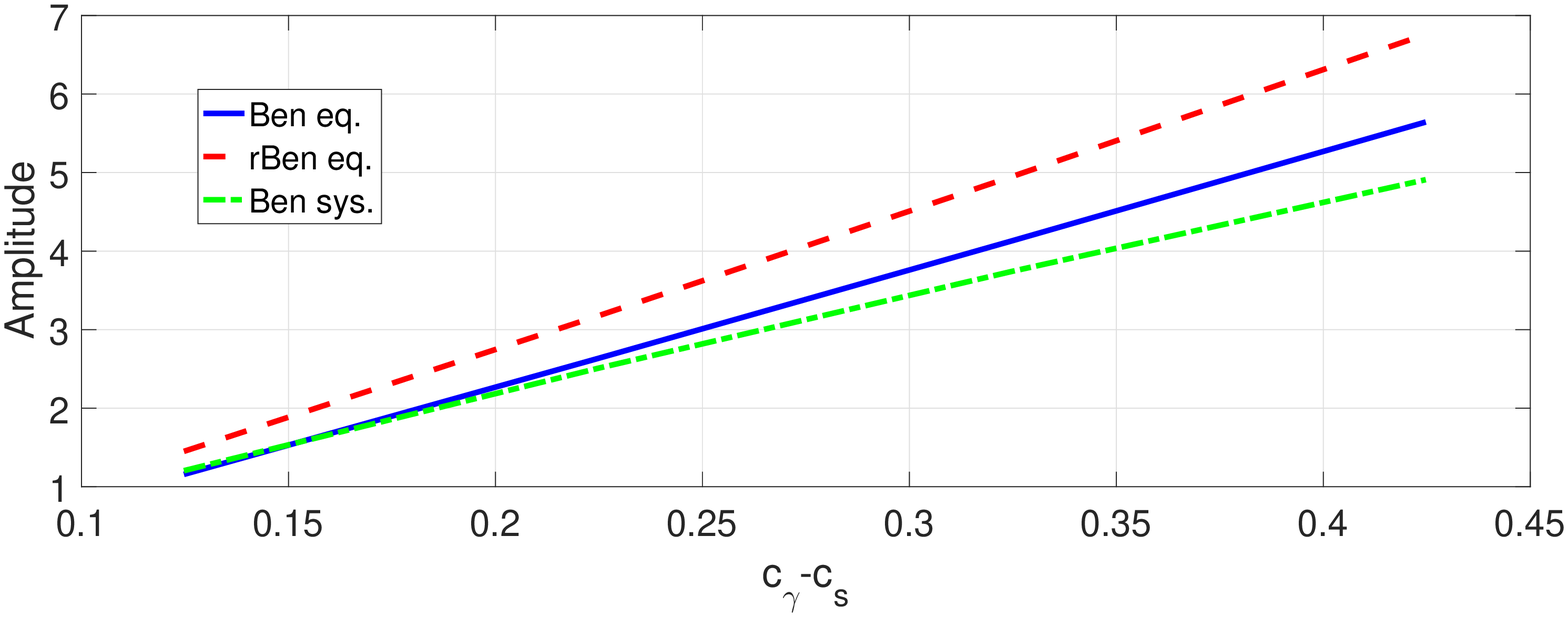}}
\subfigure[]
{\includegraphics[width=\columnwidth]{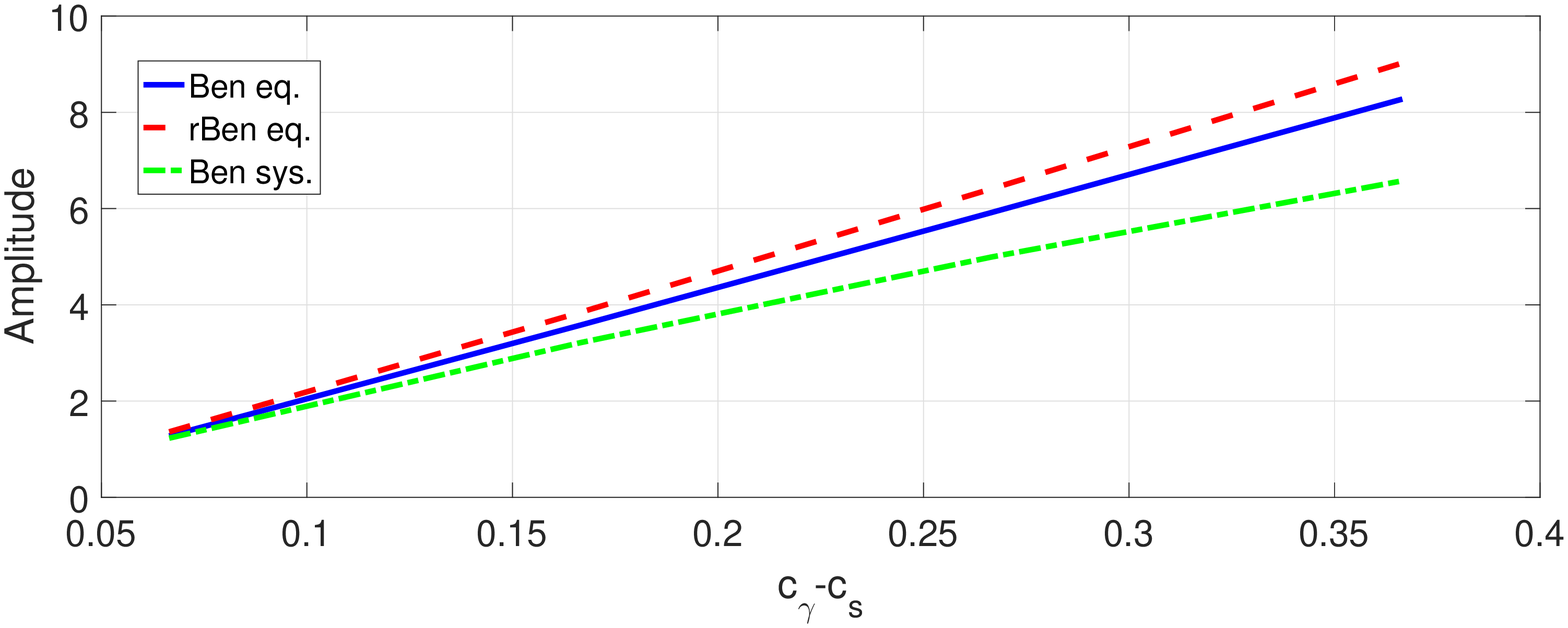}}
\subfigure[]
{\includegraphics[width=\columnwidth]{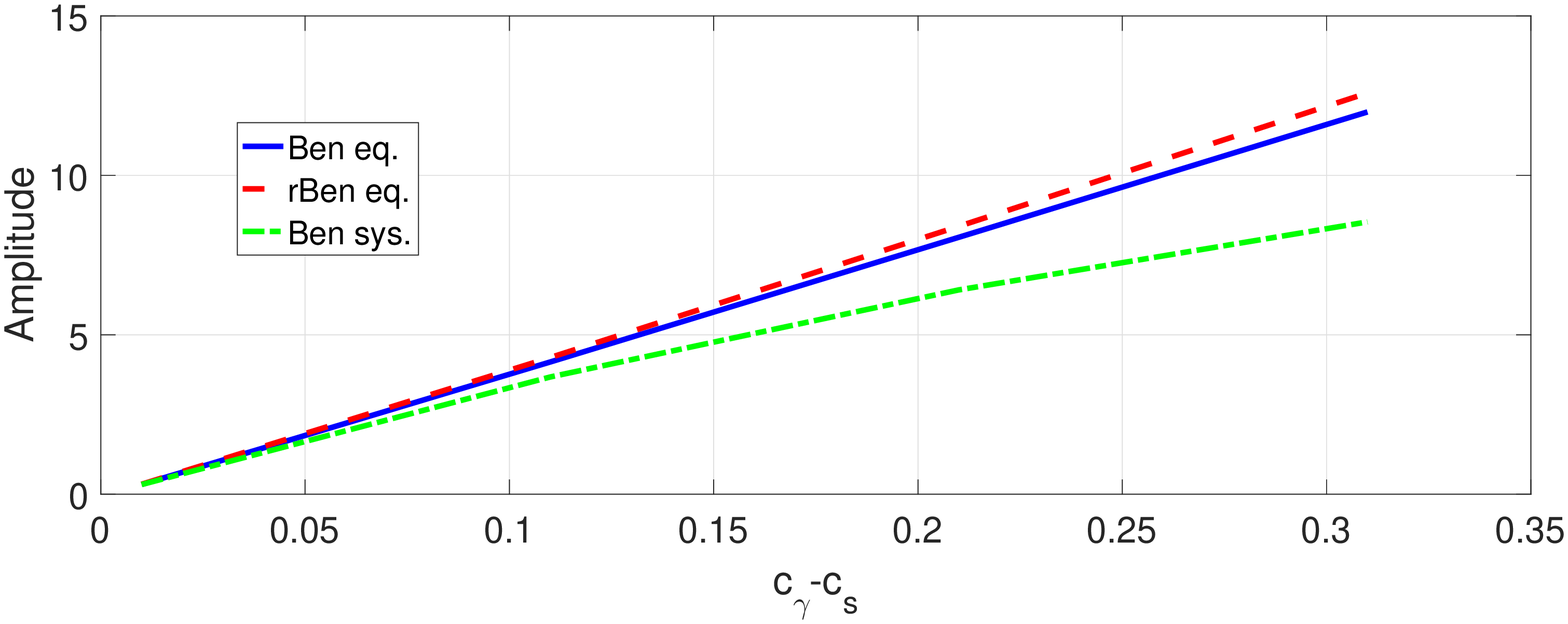}}
\caption{Speed-amplitude relation. (a) $\gamma=0.4$; (b) $\gamma=0.6$; (c) $\gamma=0.8$.}
\label{fig:bens_fig5}
\end{figure}

The main influence of the parameter $\gamma$ seems to be then through the apparently limiting speed $c_{\gamma}$. In Figure \ref{fig:bens_fig3} a comparison of the approximate $\zeta$ profiles of the Benjamin system, the rBenjamin equation and the Benjamin equation is made for different values of $\gamma=0.4, 0.6, 0.8$, for which $c_{\gamma}=1.2247, 0.8165, 0.5$ respectively. Note that as $\gamma$ grows, the profiles of the three models are closer, and the amplitude decreases. This is also observed in Figure \ref{fig:bens_fig4}, which displays the corresponding phase portraits. Here we can also notice the oscillatory decay of the waves.

\begin{figure}[!htbp] 
\centering
{\includegraphics[width=\columnwidth]{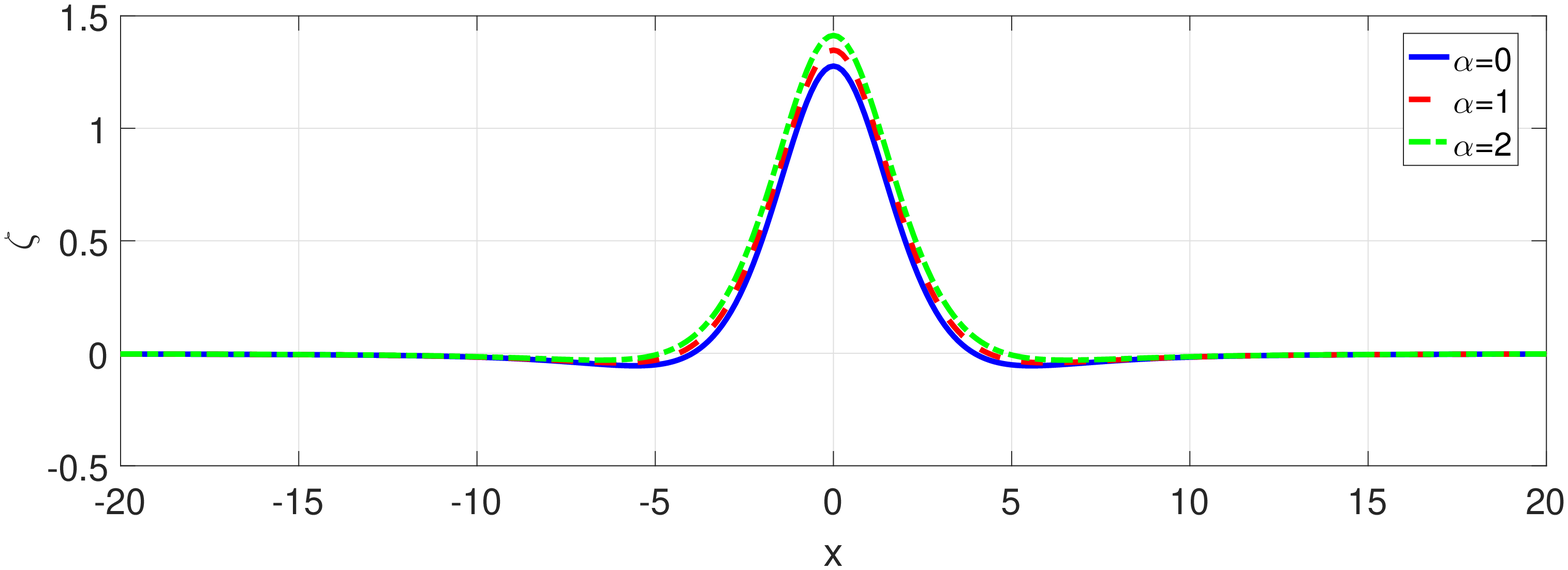}}
\caption{Approximate $\zeta$ solitary wave profiles of the rBenjamin equation for several values of $\alpha$ and $\gamma=0.6, c_s=0.75$.
}
\label{fig:bens_fig6}
\end{figure}

As mentioned before, for a fixed value of $\gamma$, the amplitude of the waves is an increasing function of $c_{\gamma}-|c_{s}|$. This is illustrated in Figure \ref{fig:bens_fig5}, which corresponds to $\gamma=0.4, 0.6$ and $0.8$. In all the cases, the model providing the largest amplitudes is the rBenjamin equation but, as $\gamma\uparrow 1$ and for larger values of $c_{\gamma}-|c_{s}|$, the behaviour of the Benjamin equation and rBenjamin equation seems to approach while the amplitudes of the profiles of the Benjamin system tend to separate from those of the corresponding for the Benjamin equation. The results in Figures \ref{fig:bens_fig3} and \ref{fig:bens_fig4} correspond to taking $\alpha=1.2$ both in (\ref{bens229}) and in (\ref{bens224}), (\ref{bens225}). Figure \ref{fig:bens_fig6} displays the profiles of the rBenjamin equation for different values of $\alpha$ and $\gamma=0.6, c_{s}=0.75$. Observe that the amplitude of the profiles is an increasing function of $\alpha$.

The similarities in the solitary waves of the three models can also be studied as follows. We generate an approximate solitary wave solution of the Benjamin equation. The profile is now considered as initial condition for two numerical methods that approximate the evolution of (\ref{bens229}) and (\ref{bens224}), (\ref{bens225}) respectively. In the case of the Benjamin system, the initial condition for the second component $u$ is given by (\ref{bens227}), (\ref{bens228}), where $\zeta$ would denote the computed solitary wave of the Benjamin equation. Then the evolution of the corresponding numerical approximation is monitored. The numerical schemes used for the simulations consist of the approximation of the corresponding periodic initial-value problem on a long enough interval with Fourier collocation discretization in space and a fourth-order, singly diagonally Runge-Kutta composition method as time integrator. Both numerical strategies were shown to have a good performance in related problems, \cite{FrutosS1992,DDM2015,DDM2019}.

\begin{figure}[htbp]
\centering
\subfigure[]
{\includegraphics[width=12cm]{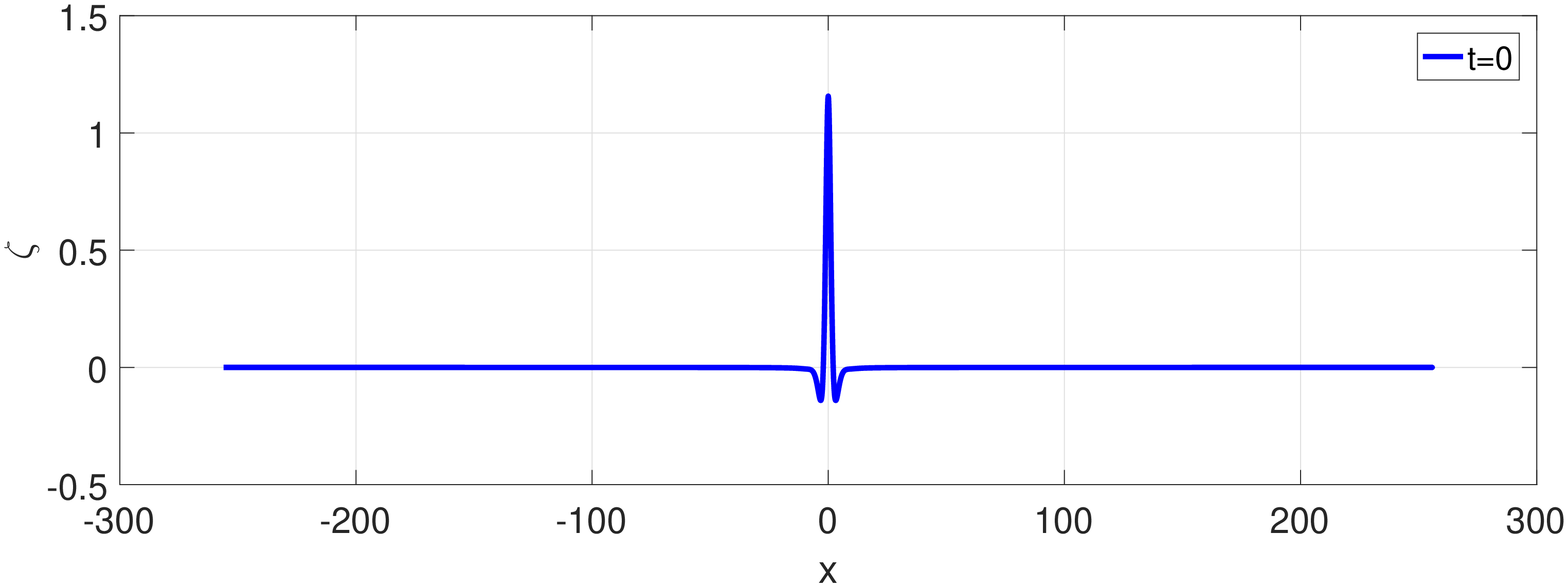}}
\subfigure[]
{\includegraphics[width=12cm]{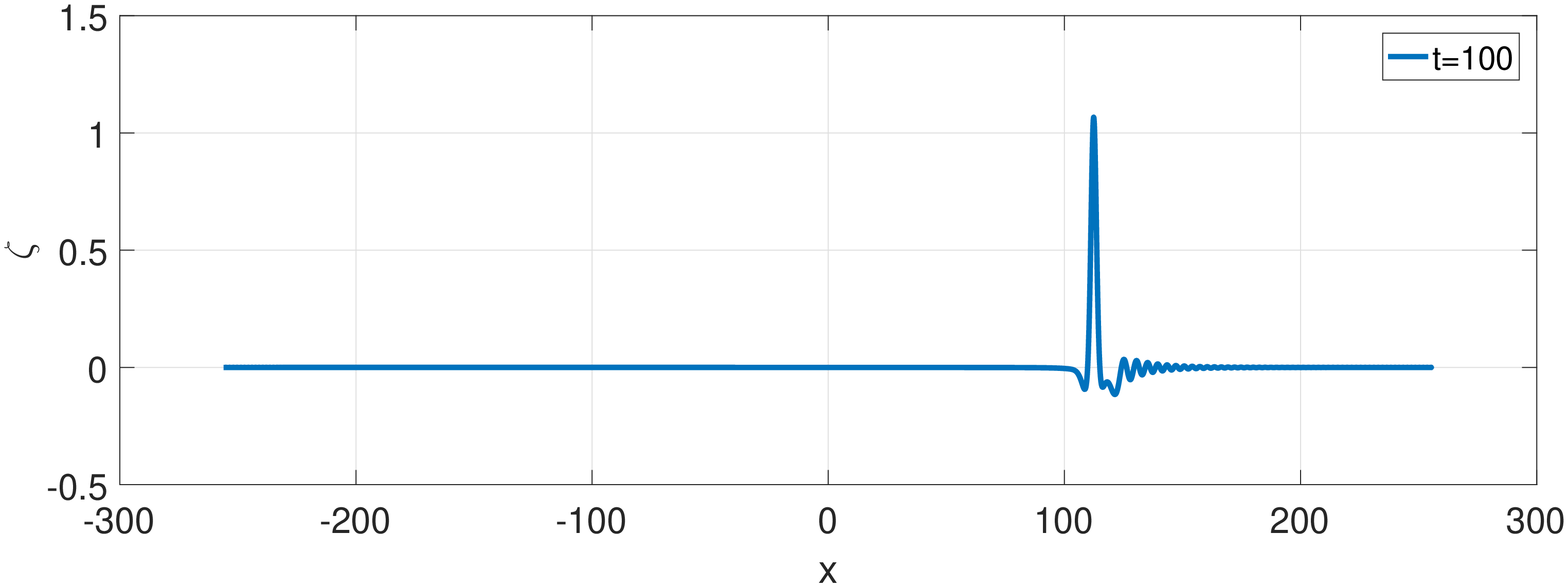}}
\subfigure[]
{\includegraphics[width=12cm]{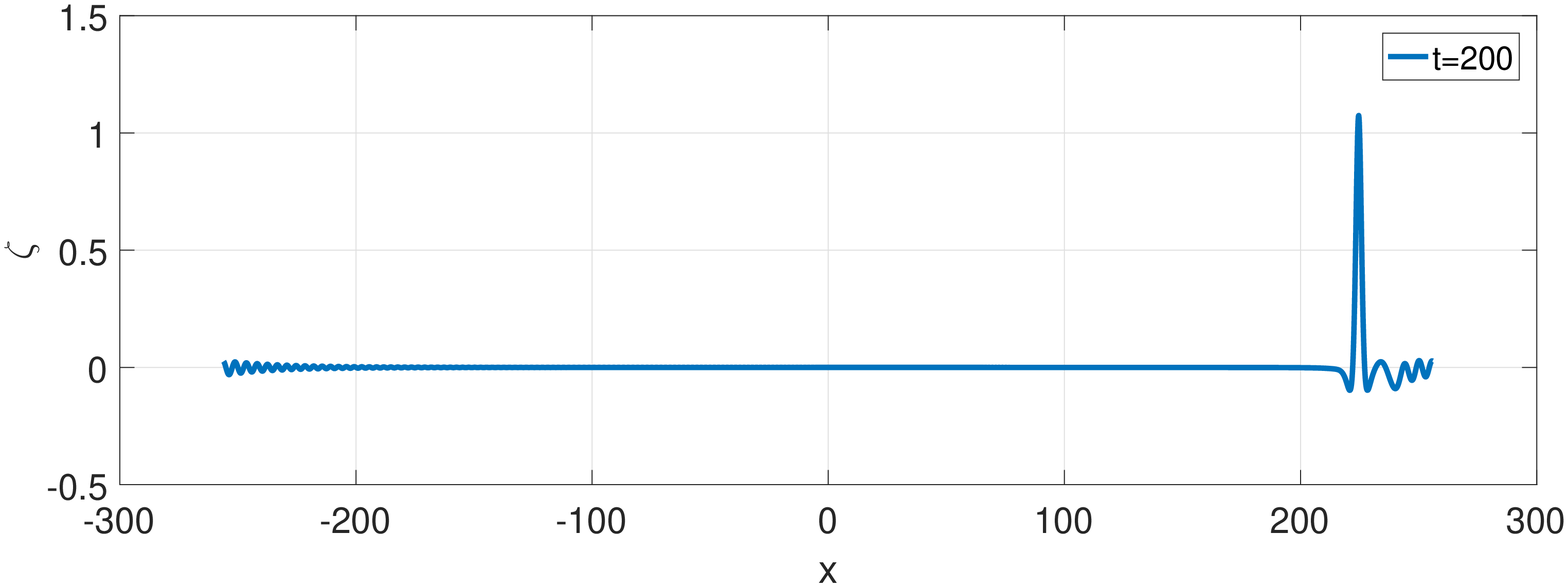}}
\subfigure[]
{\includegraphics[width=12cm]{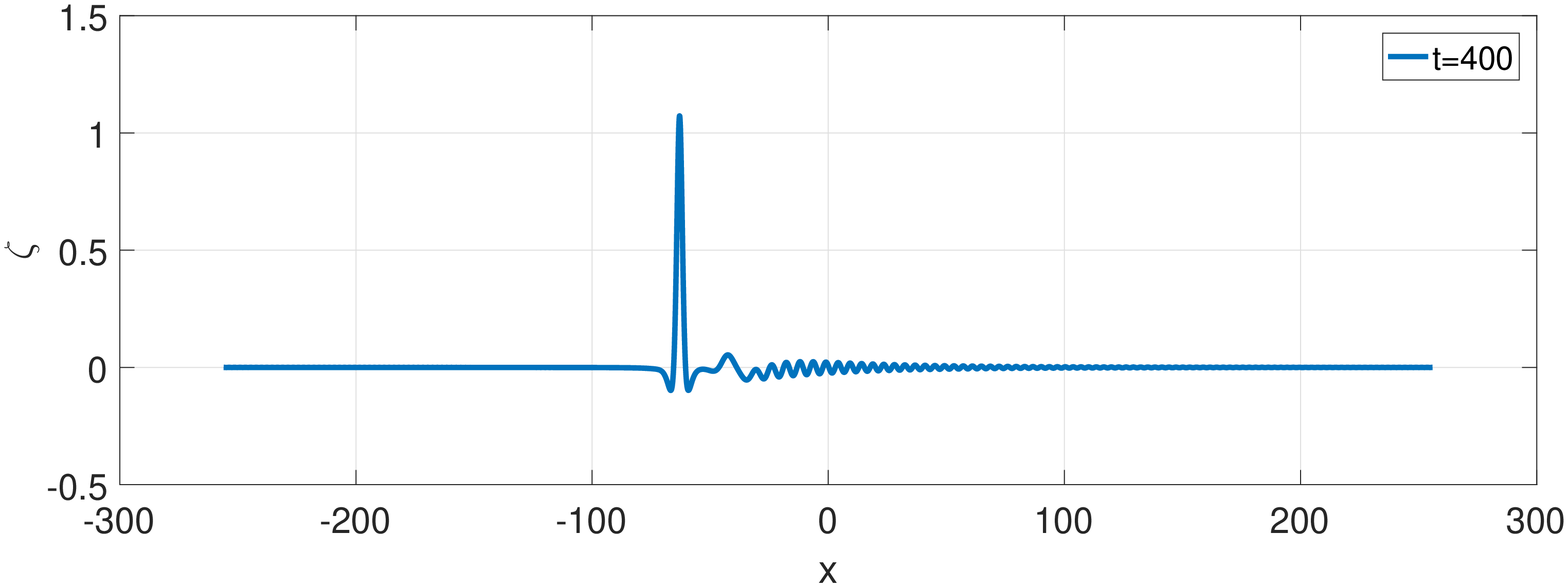}}
\caption{Numerical approximation to the rBenjamin equation from a solitary wave of the Benjamin equation as initial condition. $\gamma=0.4, c_s=1.1, \alpha=1.2$.}
\label{fig:bens_fig7}
\end{figure}

\begin{figure}[htbp]
\centering
\subfigure[]
{\includegraphics[width=12cm]{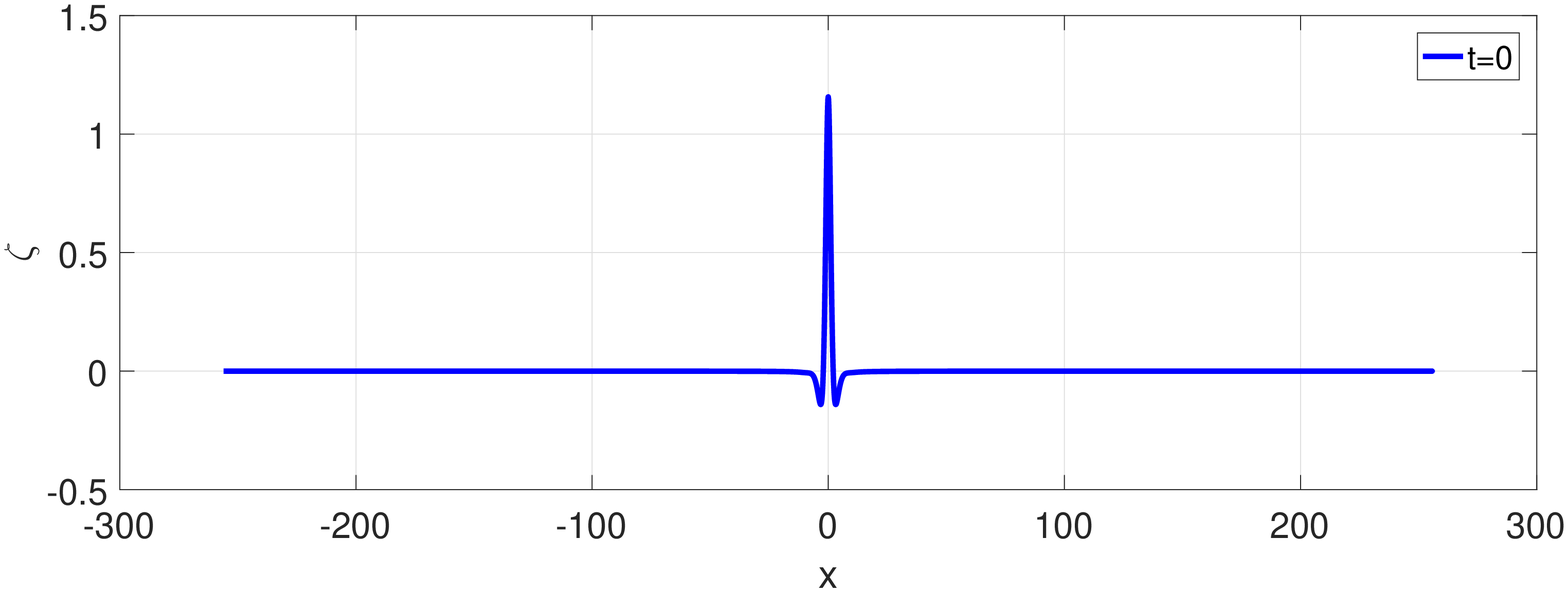}}
\subfigure[]
{\includegraphics[width=12cm]{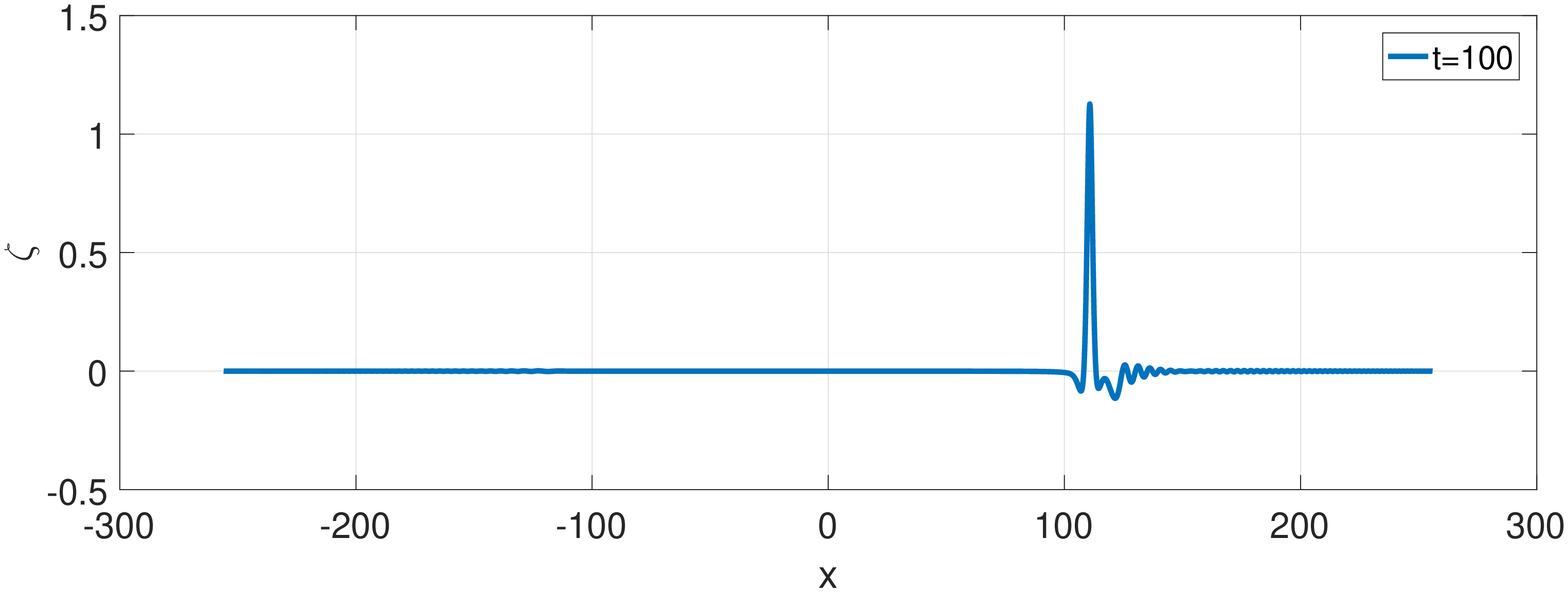}}
\subfigure[]
{\includegraphics[width=12cm]{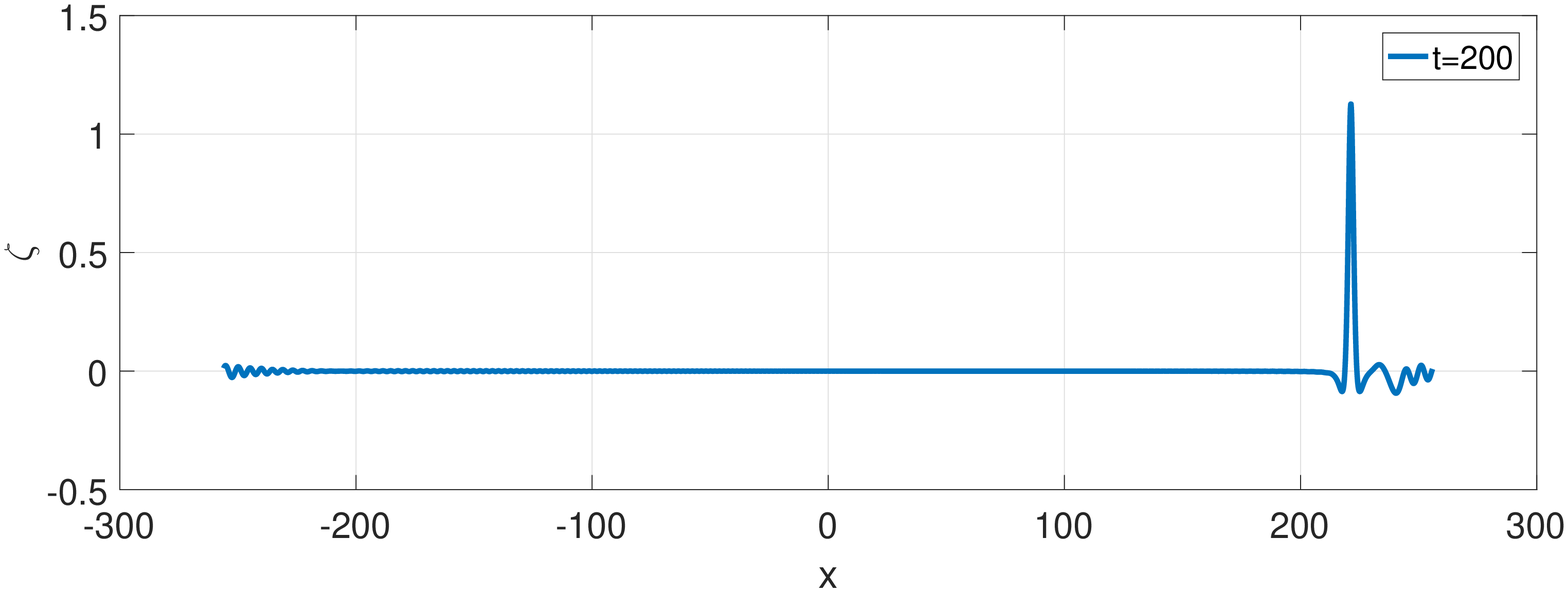}}
\subfigure[]
{\includegraphics[width=12cm]{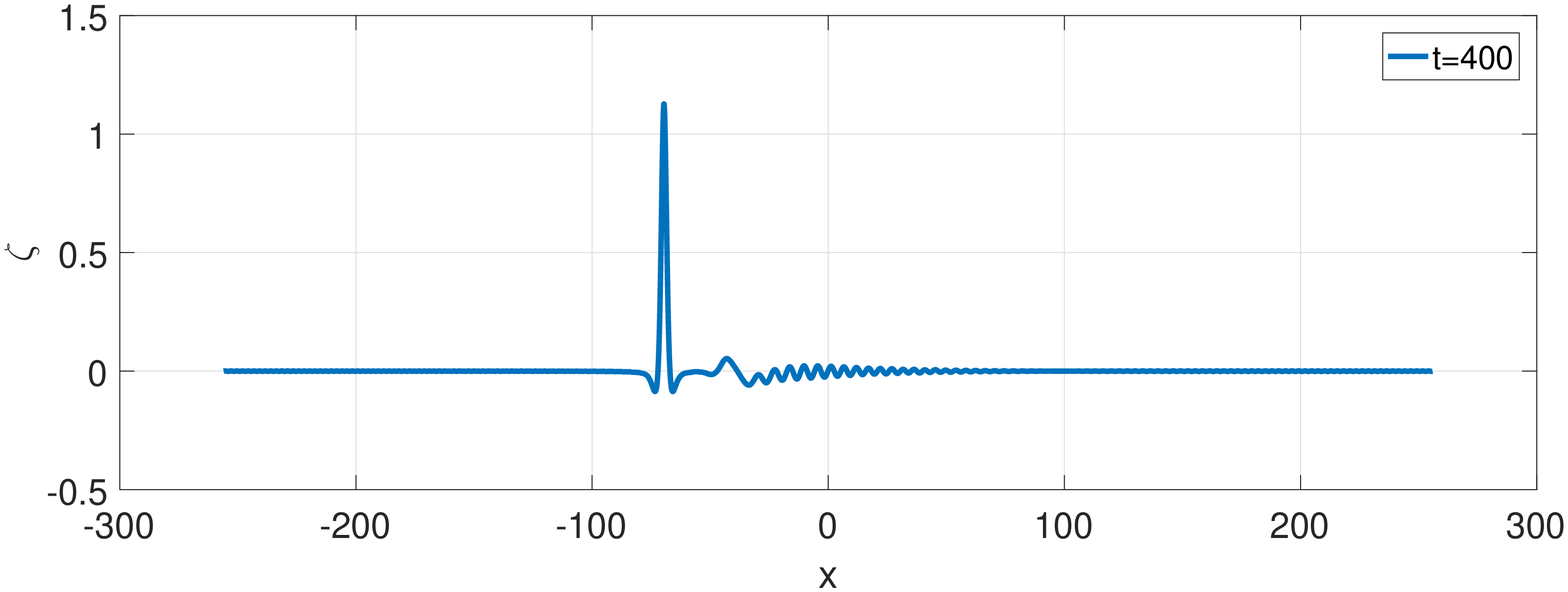}}
\caption{Numerical approximation to the Benjamin system from a solitary wave of the Benjamin equation as initial condition. $\gamma=0.4, c_s=1.1, \alpha=1.2$.}
\label{fig:bens_fig8}
\end{figure}

Taking $\gamma=0.4$ and $c_{s}=1.1$, this evolution is illustrated in Figure \ref{fig:bens_fig7} (for the rBenjamin equation with $\alpha=1.2$) and \ref{fig:bens_fig8} (for the Benjamin system with $\alpha=1.2$). 

\begin{figure}[!htbp] 
\centering
{\includegraphics[width=\columnwidth]{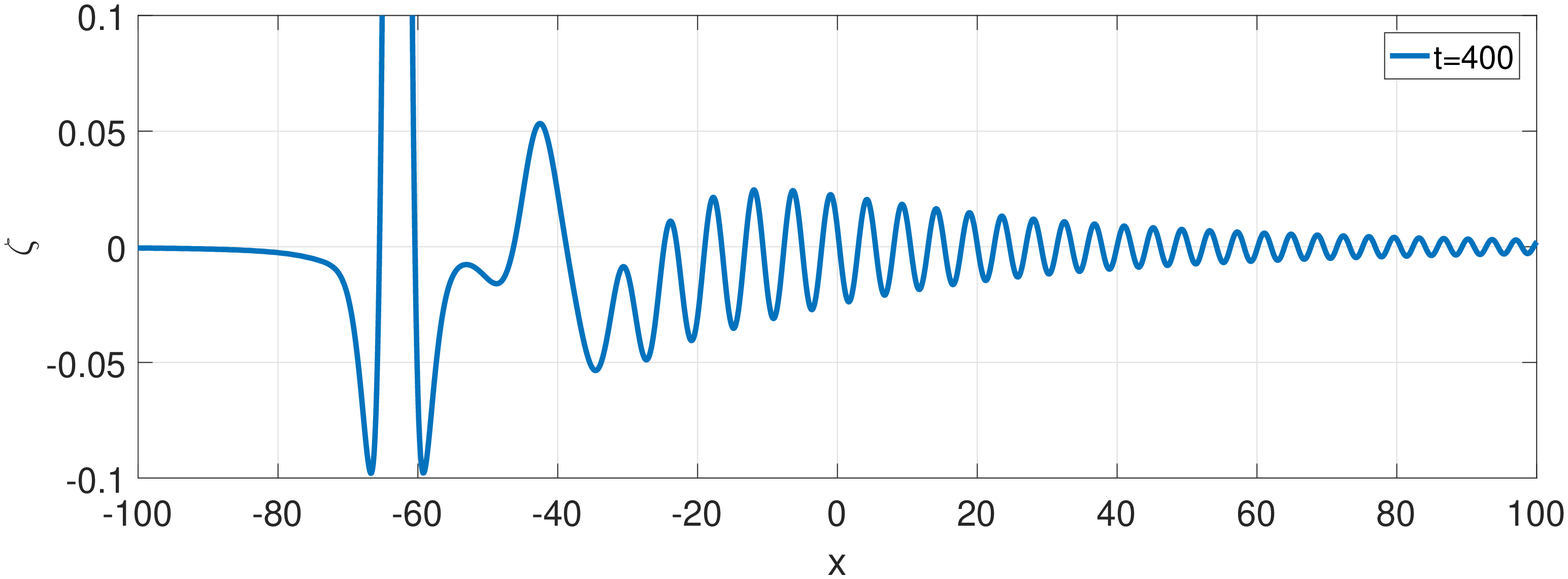}}
\caption{rBenjamin equation. Magnification of Figure \ref{fig:bens_fig7}(d)
}
\label{fig:bens_fig7m}
\end{figure}
\begin{figure}[!htbp] 
\centering
{\includegraphics[width=\columnwidth]{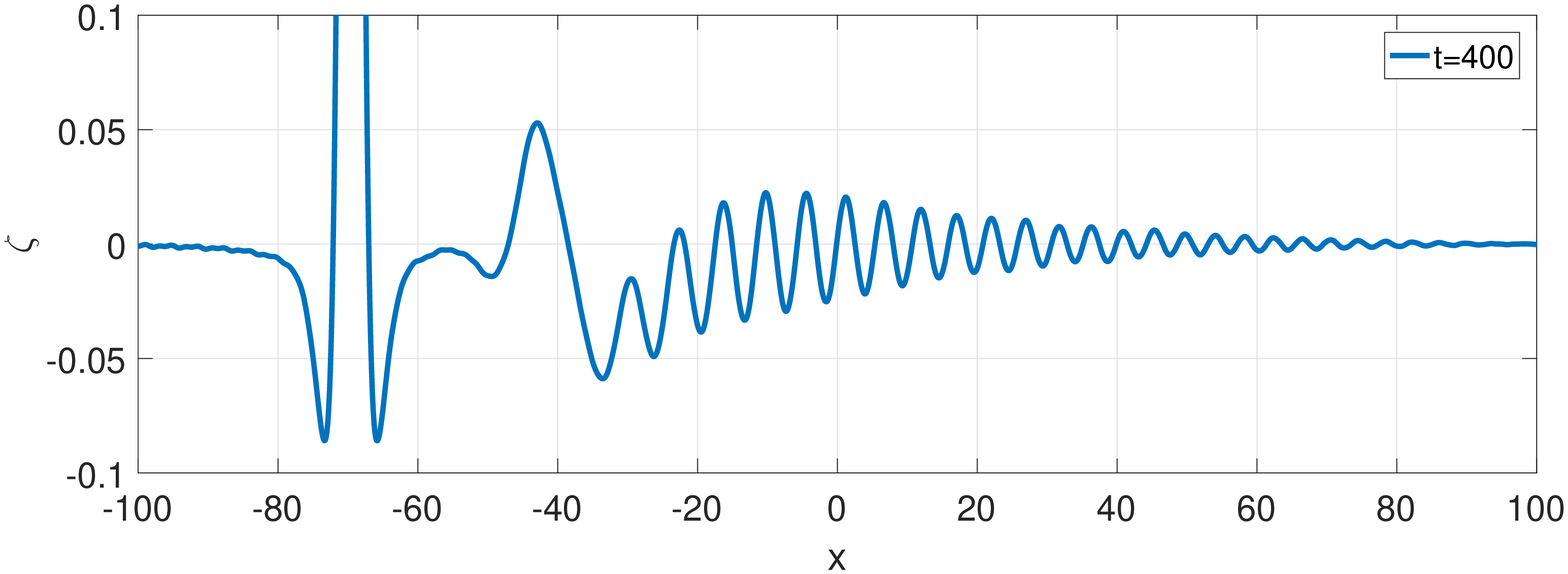}}
\caption{Benjamin system. Magnification of Figure \ref{fig:bens_fig8}(d)
}
\label{fig:bens_fig8m}
\end{figure}

The two models show a similar qualitative behaviour. The initial condition evolves into an approximate solitary wave solution of the corresponding equations along with a dispersive tail traveling in front of this main wave. (There is also a much smaller tail trailing the solitary-wave profile.) Furthermore, the formation of a small solitary wave-like structure is not discarded, in a sort of resolution property. This is suggested by the magnifications in Figures \ref{fig:bens_fig7m} and \ref{fig:bens_fig8m}. In a context of stability, the experiments suggest that the initial solitary wave solution of the Benjamin equation behaves as a small perturbation of some close solitary wave solutions of the rBenjamin equation and the Benjamin system.

\begin{figure}[!htbp] 
\centering
{\includegraphics[width=\columnwidth]{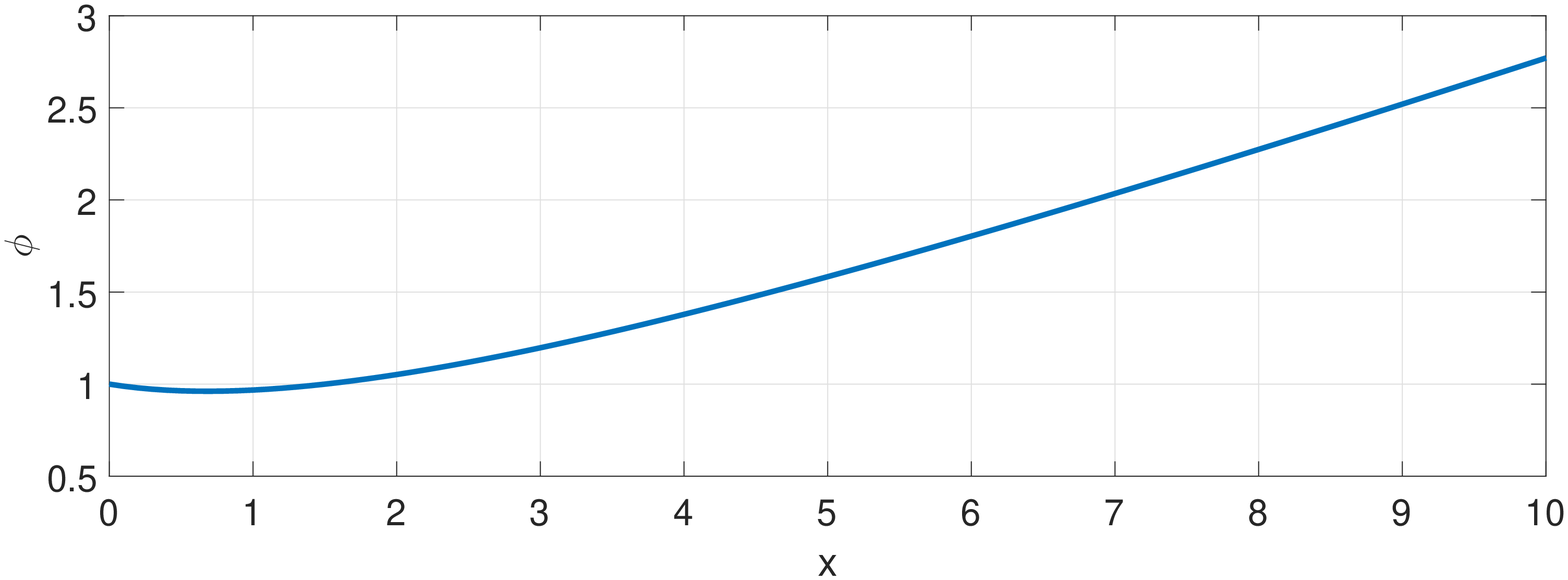}}
\caption{Function $\phi$ in (\ref{bens46}) with $\sqrt{\mu}=T=0.1$ and $\gamma=0.4, \alpha=1.2$.}
\label{fig:bens_fig9}
\end{figure}
\begin{figure}[!htbp] 
\centering
{\includegraphics[width=\columnwidth]{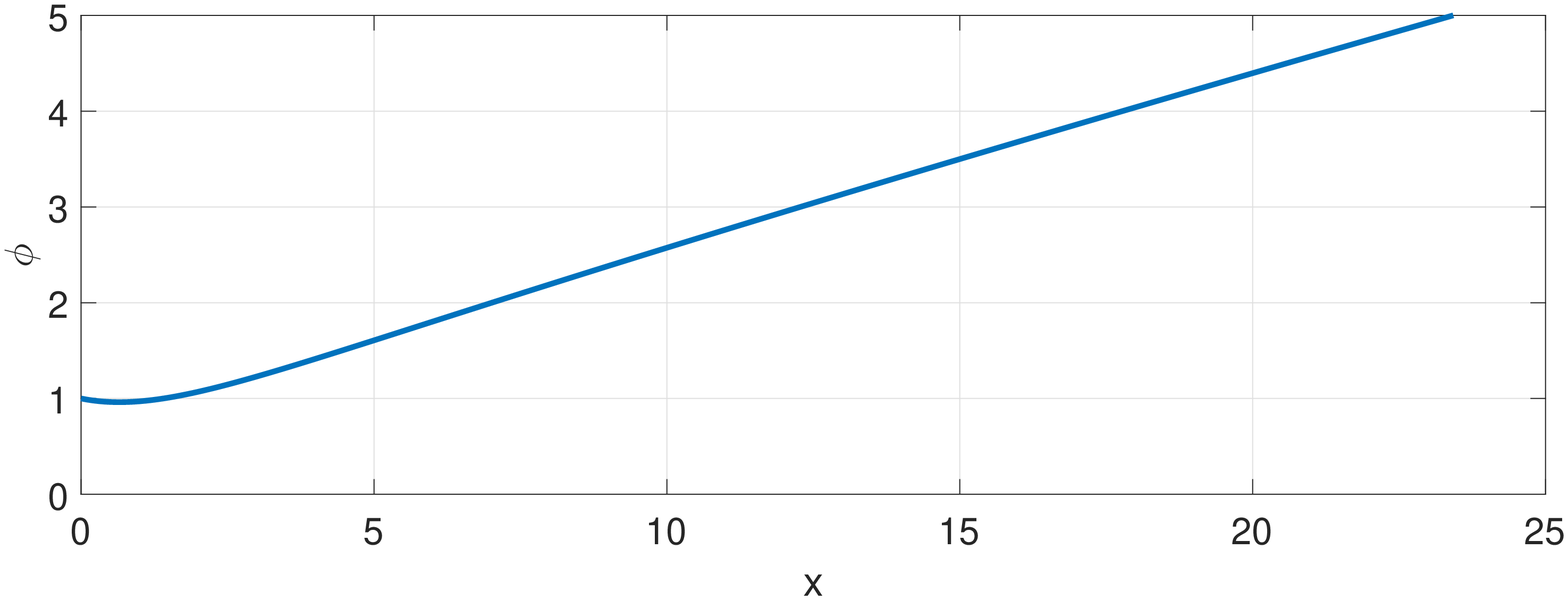}}
\caption{Function $\phi$ in (\ref{bens49}) with $\sqrt{\mu}=T=0.1$ and $\gamma=0.4, \alpha=1.2$.
}
\label{fig:bens_fig10}
\end{figure}

The structure of the dispersive tails may be studied from the corresponding linearized equations. In the case of (\ref{bens229}) and in a frame moving with the speed $c_{s}$ of the solitary wave, the equation is

\begin{eqnarray}
\left(1+\frac{\alpha\sqrt{\mu}}{\gamma}\mathcal{H}\right)(\partial_{t}-c_{s}\partial_{y})\zeta+c_{\gamma}\partial_{y}\zeta-c_{\gamma}\frac{(1-2\alpha)}{2\gamma}\sqrt{\mu}\mathcal{H}\partial_{y}\zeta&&\nonumber\\
-\frac{T}{2\sqrt{\gamma(1-\gamma)}}\partial_{y}^{3}\zeta=0,&&\label{bens45}
\end{eqnarray}
where $y=x-c_{s}t$. Plane wave solutions $\zeta(y,t)=e^{i(ky-\omega(k)t)}$ of (\ref{bens45}) will satisfy the linear dispersion relation $\omega(k)=-kc_{s}+c_{\gamma}\phi(|k|)$ where $\phi:[0,\infty)\rightarrow\mathbb{R}$ is defined as
\begin{eqnarray}
\phi(x)=\frac{1+\frac{(2\alpha-1)}{2\gamma}\sqrt{\mu}x+\frac{T}{2(1-\gamma)}x^{2}}{1+\frac{\alpha}{\gamma}\sqrt{\mu}x},\;\; x\geq 0.\label{bens46}
\end{eqnarray}
Therefore, the local phase speed relative to the speed of the solitary wave is
\begin{eqnarray*}
v(k)=\frac{\omega(k)}{k}=-c_{s}+c_{\gamma}\phi(|k|).
\end{eqnarray*}
Some properties of the function $\phi$ can explain the behaviour of the phase speed. These are collected in the following lemma.
\begin{lemma}
\label{lemm1}
The following properties of the function $\phi$ defined in (\ref{bens46}) hold:
\begin{itemize}
\item[(i)] $\phi(0)=1$ and $\displaystyle\lim_{x\rightarrow +\infty}\phi(x)=+\infty$.
\item[(ii)] $\phi$ attains a minimum at
\begin{eqnarray}
x^{*}=\frac{1}{2}\left(-b+\sqrt{b^{2}+4c}\right),\; b=\frac{2\gamma}{\alpha\sqrt{\mu}},\; c=\frac{1-\gamma}{\alpha T}, \label{bens47}
\end{eqnarray}
which satisfies $x^{*}>0$ and $\phi(x^{*})>0$ for $\mu$ small enough.
\end{itemize}
\end{lemma}

{\em Proof.} 
We write $\phi(x)=\frac{P(x)}{Q(x)}$ where
\begin{eqnarray*}
P(x)=1+\frac{(2\alpha-1)}{2\gamma}\sqrt{\mu}x+\frac{T}{2(1-\gamma)}x^{2},\;
Q(x)=1+\frac{\alpha}{\gamma}\sqrt{\mu}x,\; x\geq 0.
\end{eqnarray*}
Then elementary calculus proves (i) and the existence of $x^{*}$ given by (\ref{bens47}) where $\phi$ attains a minimum. In order to prove the last property, note that since $x^{*}>0$ then $Q(x^{*})>0$. If $2\alpha-1\geq 0$,  then $P(x^{*})>0$ and consequently $\phi(x^{*})>0$. If $2\alpha-1< 0$ we observe that $P$ attains a minimum at
\begin{eqnarray*}
x_{*}=\frac{c_{\gamma}^{2}}{2T}(1-2\alpha)\sqrt{\mu},
\end{eqnarray*}
for which, after some computations, one finds that
\begin{eqnarray*}
P(x_{*})=\frac{8\gamma^{2}T-(1-2\alpha)^{2}(1-\gamma)\mu}{8\gamma^{2}T},
\end{eqnarray*}
and since $T=O(\sqrt{\mu})$ then $P(x)\geq P(x_{*})>0$ for $\mu$ small enough and $x\geq 0$. In particular $P(x^{*})>0$ and thus $\phi(x^{*})>0$. 
$\Box$

A typical form of the function $\phi$ for the range of values of the parameters used in the numerical experiments is illustrated in Figure \ref{fig:bens_fig9}. 

Note that Lemma \ref{lemm1} implies that
\begin{eqnarray*}
v(k)>-c_{s}+c_{\gamma}\phi(x^{*}),
\end{eqnarray*}
and $v(k)>0$ for all wavenumbers $k$ if we take $c_{s}<c_{\gamma}\phi(x^{*})$. Even if this is not satisfied, we observe that, due to (i) of Lemma \ref{lemm1}, we have $\phi(|k|)>1$ from some value of $|k|$ (see Figure \ref{fig:bens_fig9}) and then
\begin{eqnarray*}
v(k)>-c_{s}+c_{\gamma}.
\end{eqnarray*}
Therefore, if $c_{s}<c_{\gamma}$ then $v(k)>0$ from some value of $|k|$; this means that most of the solution components $\zeta(y,t)=e^{i(ky-\omega(k)t)}$ is leading the solitary pulse, cf. Figure \ref{fig:bens_fig7}.

In the case of the Benjamin system (\ref{bens224}), (\ref{bens225}), the corresponding linearized equations are
\begin{eqnarray}
\left(1+\frac{\alpha\sqrt{\mu}}{\gamma}\mathcal{H}\right)(\partial_{t}-c_{s}\partial_{y})\zeta+\frac{1}{\gamma}\partial_{y}\left(1-(1-\alpha)\frac{\sqrt{\mu}}{\gamma^{2}}\mathcal{H}\right){u}&=&0,\label{bens48a}\\
\partial_{t}{u}+\partial_{y}\left((1-\gamma)-T\partial_{y}^{2}\right)\zeta&=&0,\label{bens48b}
\end{eqnarray}
System (\ref{bens48a}), (\ref{bens48b}) can be reduced to
\begin{eqnarray*}
\left(1+\frac{\alpha\sqrt{\mu}}{\gamma}\mathcal{H}\right)(\partial_{t}-c_{s}\partial_{y})^{2}\zeta-c_{\gamma}^{2}\partial_{y}^{2}\left(1-(1-\alpha)\frac{\sqrt{\mu}}{\gamma^{2}}\mathcal{H}\right)\zeta&&\\
-\frac{T}{\gamma}\partial_{y}^{4}\left(1-(1-\alpha)\frac{\sqrt{\mu}}{\gamma^{2}}\mathcal{H}\right)\zeta=0&&.
\end{eqnarray*}
Then the dispersion relation, relative to the speed $c_{s}$ of the solitary wave, has the form $\omega_{\pm}(k)=-kc_{s}\pm kc_{\gamma}\phi(|k|)$, where now
$\phi:[0,\infty)\rightarrow\mathbb{R}$ is given by
\begin{eqnarray}
\phi(x)=\left(\frac{(1+\frac{(\alpha-1)}{\gamma}\sqrt{\mu}x)(1+\frac{T}{(1-\gamma)}x^{2})}{1+\frac{\alpha}{\gamma}\sqrt{\mu}x}\right)^{1/2},\;\; x\geq 0.\label{bens49}
\end{eqnarray}
We observe (see Figure \ref{fig:bens_fig10}) that the function $\phi$ in (\ref{bens49}) also satisfies (i) of Lemma \ref{lemm1}. If $|c_{s}|<c_{\gamma}$, then from some $|k|$ it holds that
\begin{eqnarray*}
v_{+}(k)=\frac{\omega_{+}(k)}{k}=-c_{s}+c_{\gamma}\phi(|k|)>-c_{s}+c_{\gamma}>0,
\end{eqnarray*}
and
\begin{eqnarray*}
v_{-}(k)=\frac{\omega_{-}(k)}{k}=-c_{s}-c_{\gamma}\phi(|k|)<-c_{s}-c_{\gamma}<0.
\end{eqnarray*}
Thus, as before,  most of the wave components of the dispersive tail travels leading the solitary wave with speed $c_{s}$.
\section{Concluding remarks}
\label{sec:sec4}
The present paper introduces a two-dimensional asymptotic model for the propagation of internal waves in a two-layer system of fluids with rigid lid condition for the upper layer and a lower layer of infinite depth (or much larger than that of the upper layer). Furthermore, the model takes into account both gravity and capillary effects at the interface. The derivation is carried out reformulating first  the corresponding Euler equations by using the nonlocal operators considered in \cite{BLS2008}. Then asymptotic expansions of these operators, consistent with the physical regime of the Benjamin model, lead to a bi-directional system for the deviation of the interface and the velocity variables. In the one-dimensional, uni-directional case, the system can be reduced to a family of regularized Benjamin equations which contains, as particular case, a version of the model derived by Benjamin, \cite{Benjamin1967,Benjamin1992}. 

Some mathematical properties of the new models are also discussed: linear well-posedness, existence of conserved quantities and a computational study of comparison of solitary wave solutions for the one-dimensional Benjamin system, the regularized Benjamin equation and the usual Benjamin equation. The results obtained about these mathematical aspects will additionally serve us as starting point for a deeper study of the models, concerning topics of local and global well-posedness, existence of solitary waves and dynamics of the equations. 
%
%





\begin{thebibliography}{ablowitzm}
\bibitem{ABR} J. P. Albert, J. L. Bona, J. M. Restrepo, Solitary-wave solutions of the Benjamin equation, SIAM J. Appl. Math. 59 (1999) 2139-2161.
\bibitem{AlvarezD2015} J. \'Alvarez, A. Dur\'an, Petviashvili type methods for traveling wave computations: II. Acceleration with vector extrapolation methods, Math. Comput. Simul., 123 (2016) 19-36.
\bibitem{AnguloS2020} J. Angulo-Pava, J.-C. Saut, Existence of solitary wave solutions for internal waves in two-layer systems, Quart. Appl. Math. 78 (2020), 75-105.
\bibitem{Benjamin1967} T. B. Benjamin, Internal waves of permanent form in fluids of great depth, J. Fluid Mech. 29 (1967) 559-592.
\bibitem{Benjamin1992} T.B. Benjamin, A new kind of solitary wave, J. Fluid Mech. 245 (1992) 401-411.
\bibitem{Benjamin1996} T. B. Benjamin, Solitary and periodic waves of a new kind, Philos. Trans. Roy. Soc. London Ser. A 354 (1996) 1775-1806.
\bibitem{BenjaminB} T. B. Benjamin, T. J. Bridges, Reappraisal of the Kelvin-Helmholtz problem. Part 1. Hamiltonian structure, J. Fluid Mech. 333 (1997) 301-325.
\bibitem{Bona1981} J. L. Bona, On solitary waves and their role in the evolution of long waves, In {\em Applications of Nonlinear Analysis in the Physical Sciences} (ed. H. Amann, N. Bazlev, K. Kirchg\"{a}ssner) Pitman, London, pp. 183-205, 1981.
\bibitem{BonaDM2020} J. L. Bona, A. Dur\'an, D. Mitsotakis, Solitary-wave solutions of Benjamin-Ono and other systems for internal waves. I.  Approximations, 2020, in Press.
\bibitem{BLS2008}
J. L. Bona, D. Lannes, J. C. Saut, Asymptotic models for internal waves, J. Math. Pures Appl., 89 (2008), 538-566.
\bibitem{CA} D. C. Calvo, T. R. Akylas, On interfacial gravity-capillary solitary waves of the Benjamin type and their stability, Phys. Fluids 15 (2003) 1261-1270.
\bibitem{CraigS} W. Craig, C. Sulem, Numerical simulation of gravity waves, J. Comput. Phys. 108 (1993) 73-83.
\bibitem{DiasI1996} F. Dias, G. Iooss, Capillary-gravity interfacial waves in infinite depth, Eur. J. Mech. B Fluids 15 (1996) 367-393.
\bibitem{DiasMV1996} F. Dias, D. Menasce, J. M. Vanden-Broeck, Numerical study of capillary-gravity solitary waves, Eur. J. Mech. B Fluids 15(1) 17-36.
\bibitem{DDM2015} V. A. Dougalis, A. Dur\'an, D. E. Mitsotakis, Numerical solution of the Benjamin equation, Wave Motion 52 (2015) 194-215.
\bibitem{DDM2019} V. A. Dougalis, A. Dur\'an, D. E. Mitsotakis, Numerical approximation to Benjamin type equations. Generation and stability of solitary waves, Wave Motion 85 (2019) 34-56. 
\bibitem{DougalisMS2007} V. A. Dougalis, D. E. Mitsotakis, J.-C. Saut, On some Boussinesq systems in two space dimensions: Theory and numerical analysis, ESAIM: Math. Modelling and Numer. Anal. 41 (2007) 825-854.
\bibitem{FrutosS1992} J. de Frutos, J. M. Sanz-Serna, An easily implementable fourth-order method for the time integration of wave problems, J. Comput. Phys., 103 (1992) 160-168.
\bibitem{HelfrichM} K. R. Helfrich, W. K. Melville, Long nonlinear internal waves, {Annual Review of Fluid Mechanics} 38 (2006) 395-425.
\bibitem{KB} H. Kalisch, J. L. Bona, Models for internal waves in deep water, Discret. Contin. Dyn. Syst. 6 (2000) 1-22.
\bibitem{Kalisch2007} H. Kalisch, Derivation and comparison of model equations for interfacial capillary-gravity waves in deep water, Math. Comput. Simul. 74 (2007) 168-178.
\bibitem{KenigPV1991} C.E. Kenig, G. Ponce, L. Vega, Well-posedness of the initial value problem for the Korteweg-de Vries equation, J. Amer. Math. Soc. 4 (1991) 323-347.
\bibitem{KenigPV1993} C.E. Kenig, G. Ponce, L. Vega, Well-posedness and scattering results for the generalized Korteweg-de Vries equation via contraction principle, Comm. Pure Appl. Math. 46 (1993) 527-620.
\bibitem{KimA2005} B. Kim, T. R. Akylas, On gravity-capillary lumps, J. Fluid Mech. 540 (2005) 337-351.
\bibitem{KimA2006} B. Kim, T. R. Akylas, On gravity-capillary lumps. Part 2. Two-dimensional Benjamin equation, J. Fluid Mech. 557 (2006) 237-256.
\bibitem{KoopB1981} C. G. Koop, G. Butler, An investigation of internal solitary waves in a two-fluid system, J. Fluid Mech. 112 (1981) 225-251.
\bibitem{LagetD1997} O. Laget, F. Dias, Numerical computation of capillary-gravity interfacial solitary waves, J. Fluid Mech. 349 (1997) 221-251.
\bibitem{Lannes} D. Lannes, A stability criterion for two-fluid interfaces and applications, Arch. Ration. Mech. Anal. 208 (2013) 481-567.  
\bibitem{Linares1999} F. Linares, $L^{2}$ global well-posedness of the initial value problem associated to the Benjamin equation, J. Diff. Eq., 152 (1999) 377-393.
\bibitem{LinaresS2005} F. Linares, M. Scialom, On generalized Benjamin type equations, Disc. Cont. Dyn. Sys., 12(1) (2005) 161-174.
\bibitem{pelinovskys} {D. E. Pelinovsky and Y. A.
Stepanyants}, { Convergence of Petviashvili's iteration method for
numerical approximation of stationary solutions of nonlinear wave
equations}, { SIAM J.~Numer.\ Anal.} { 42} (2004) 1110-1127.
\bibitem{Petv1976} { V. I. Petviashvili} { Equation of
an extraordinary soliton}, { Soviet J.~Plasma Phys.} { 2} (1976)
257-258.
\bibitem{sidi} A. Sidi, Vector Extrapolation Methods with Applications, SIAM Philadelphia, 2017.
\bibitem{sidifs} A. Sidi, W. F. Ford, D. A. Smith, Acceleration of convergence of vector sequences, SIAM J. Numer. Anal., 23 (1986) 178-196.
\bibitem{smithfs} D. A. Smith, W. F. Ford, A. Sidi, Extrapolation methods for vector sequences, SIAM Rev., 29 (1987) 199-233.
\bibitem{Xu2012} L. Xu, Intermediate long wave systemsa for internal waves, Nonlinearity 25 (2012) 597-640.
\bibitem{Whitham} G. B. Whitham, Linear and Nonlinear Waves, Wiley, New York, 1974.
\bibitem{Zakharov} V. E. Zakharov, stability of periodic waves of finite amplitude on the surface of a deep fluid, J. Applied Mechanics and technical Physics 9 (1968) 190-194.
\end{thebibliography}
\end{document}